\documentclass[journal]{IEEEtran}
\usepackage{amsmath}
\usepackage{amssymb}
\usepackage{makecell}
\usepackage{amsthm} 

\usepackage{cite}

\newcommand{\m}{\boldsymbol}

\newcommand{\mbb}[1]{\mathbb{#1}}

\newtheorem{thm}{Theorem}
\newtheorem{lem}{Lemma}

\newtheorem{remark}{Remark}

\newtheorem{assumption}{Assumption}

\DeclareMathOperator{\diag}{diag}

\usepackage{mathtools}
\usepackage{multirow}
\usepackage{multicol}
\usepackage[linesnumbered,lined,boxed,commentsnumbered,ruled,longend]{algorithm2e}
\DeclarePairedDelimiter{\norm}{\lVert}{\rVert}

\usepackage{amsthm}
\usepackage{tikz}
\usetikzlibrary{matrix,positioning,decorations.pathreplacing}
\usepackage{array}
\usepackage{filecontents}
\usepackage{cases}
\usepackage[export]{adjustbox}

\usepackage{index}     
\usepackage[linktoc=all]{hyperref}
\hypersetup{
	colorlinks = true,
	citecolor=black,
	filecolor=black,
	linkcolor=blue,
	urlcolor=black
}

\usepackage{lineno}
\begin{document}
	
	%
	%
	
	
	\title{A New Derivative-Free Linear Approximation for Solving the Network Water Flow Problem with Convergence Guarantees}
	\author{Shen Wan$\text{g}^\dagger$, Ahmad F. Tah$\text{a}^\dagger$, Lina Sel$\text{a}^\ddagger$, Marcio H. Giacomon$\text{i}^\ast$, Nikolaos Gatsi$\text{s}^\dagger$ \thanks{
		$^\dagger$Department of Electrical and Computer Engineering, The University of Texas at San Antonio, TX 78249, $^\ddagger$Department of Civil, Architectural and Environmental Engineering, Cockrell School of Engineering, The University of Texas at Austin, $^\ast$Department of Civil and Environmental Engineering, The University of Texas at San Antonio. Emails: mvy292@my.utsa.edu, ahmad.taha@utsa.edu, linasela@utexas.edu, \{marcio.giacomoni, nikolaos.gatsis\}@utsa.edu. This material is based upon work supported by the National Science Foundation under Grants CMMI-DCSD-1728629 and 1847125.  This work was also supported by the University of Texas at Austin Startup Grant and by Cooperative Agreement No. 83595001 awarded by the U.S. Environmental Protection Agency (EPA) to The University of Texas at Austin. }}
	\maketitle
	\begin{abstract}
		Addressing challenges in urban water infrastructure systems including aging infrastructure, supply uncertainty, extreme events, and security threats, depend highly on water distribution networks modeling emphasizing the importance of realistic assumptions, modeling complexities, and scalable solutions. In this study, we propose a derivative-free, linear approximation for solving the network water flow problem (WFP). The proposed approach takes advantage of the special form of the nonlinear head loss equations and, after the transformation of variables and constraints, the WFP reduces to a linear optimization problem that can be efficiently solved by modern linear solvers. Ultimately, the proposed approach amounts to solving a series of linear optimization problems. We demonstrate the proposed approach through several case studies and show that the approach can model arbitrary network topologies and various types of valves and pumps, thus providing modeling flexibility. Under mild conditions, we show that the proposed linear approximation converges. We provide sensitivity analysis and discuss in detail the current limitations of our approach and suggest solutions to overcome these. All the codes, tested networks, and results are freely available on Github for research reproducibility.
	\end{abstract}
	

	%
	%
	
	%
	
	
	%
	%
	%
	%

\section{Introduction}~\label{sec:literature}
Water distribution networks (WDNs) are complex, large-scale critical infrastructure responsible for providing safe drinking water to the continuously growing population. In the U.S., public water utilities serve approximately 95\% of the total population \cite{EPA_2017}. WDNs are composed of numerous elements such as pipes, valves, tanks, and pumps that transport water from a few points of water supply to numerous water consumers. Aging infrastructure, supply uncertainty, growing population, extreme events, and security threats, pose mounting challenges on urban water infrastructure \cite{black_veatch,AWWA_2017}. Addressing these challenges depend highly on WDN modeling and the validity of these models. Furthermore, taking advantage of technological advances and integrating smart sensing and actuation with physically-based models for operations and management of urban water systems heavily relies on the WDN models.
Modeling WDNs involves solving the water flow problem (WFP), which is governed by the \textit{linear} flow continuity and \textit{nonlinear} energy conservation \cite{todini2013}, and finding the \textit{flow} through each pump, valve, and pipe and \textit{head} at each node given network characteristics, status of pumps and valves, initial head of tanks and reservoirs, and demand at each node of the WDN. Notably, realistic assumptions, modeling complexities, and the inherently large-scale of WDNs, emphasize the importance of scalable water flow modeling solutions in the context of managing modern WDNs.

The literature of solving the WFP as well as other related problem formulations is rich and briefly summarized next. {The main approaches for solving the WFP are based on Hardy-Cross ~\cite{cross1936analysis}, Newton-Raphson ~\cite{martin1963application,liu1969numerical,epp1970efficient,wood1981reliability,todini1987s}}, linearization~\cite{wood1972hydraulic,jeppson1976analysis,isaacs1980linear,wood1993hydraulic,price2012iterative,moosavian2017multilinear}, optimization~\cite{arora1976flows,collins1978solving},  gradient-based~\cite{todini1987s,todini2013}, graph decomposition \cite{ulanicki,deuerlein2008,deuerlein2009,DIAO2014253} and more recently,  fixed-point methods~\cite{Zhang2017a,Hafez-FixedPointWDSA}.
These methods can be classified as primarily relying on iterative updating, decomposition methods, or optimization-based formulations, and differ in terms of their modeling limitations and complexity, handling non-linearities, and convergence speed, as discussed in the following paragraphs.

The first classical approach is attributed to \cite{cross1936analysis} that developed a loop-based method for solving the WFP suitable for small networks and hand-calculations. 
\cite{martin1963application} first applied the Newton-Raphson method by modeling all the equations in terms of nodal heads and obtaining the solution via successive iterations. Slow convergence and large oscillations during iterations are the two main disadvantages of the proposed approach. 
Later, \cite{liu1969numerical} proposed a simplified version of the Newton-Raphson method via decomposition of the Jacobian matrix into diagonal and non-diagonal matrices, which simplified the solution of the equations. However, the method suffers from convergence issues if the initial guess is not carefully chosen. 
A linearization method was proposed by \cite{wood1972hydraulic} in terms of link flow equations, where the nonlinear energy equations were linearized and updated in each iteration \cite{wood1972hydraulic}. \cite{wood1993hydraulic} later extended the link flow model using extended Taylor series. 
\cite{jeppson1976analysis} reformulated the nonlinear energy equations for each loop in terms of flow adjustment factors and proposed a linearization method using the standard Taylor series expansion, which was then solved iteratively using the Newton Raphson method.
\cite{isaacs1980linear} proposed a linearization method based on nodal heads providing a simpler model and symmetry of coefficient matrix compared with \cite{wood1972hydraulic}. 
Notably, the global gradient algorithm ~\cite{todini1988gradient} implemented in the EPANET software~\cite{rossman2000epanet}, which utilizes the Newton-Raphson solution approach for solving the nonlinear system of equations, is the most widely used method for solving the WFP \cite{burger2016}. 
\cite{giustolisi2011computationally} proposed an enhanced global gradient method to accelerate the convergence process for large-scale networks while preserving the accuracy of the solution. \cite{moosavian2017multilinear} derived a multi-linear method to improve the convergence rate of ~\cite{wood1972hydraulic} and ~\cite{todini1988gradient}, where the nonlinear energy equations are linearized based on the maximum and minimum allowable flow rate in pipes and the solution is iteratively updated in the successive iterations. To further accelerate and improve convergence several recent works have proposed intricate algorithms to exploit network structure in the computational procedure including careful selection of network loops~\cite{alvarruiz2015improving, deuerlein2015fast, vasilic2018improved}, selection and decomposition of network trees and forest~\cite{simpson2012forest,elhay2014reformulated}. Recently, the uniqueness of WFP is discussed in~\cite{singh2019flow}.

An alternative approach for solving the WFP is by formulating the problem as nonlinear but convex optimization problem, i.e. the content problem which is constrained by linear mass balance equations minimizing network content or the unconstrained dual problem minimizing the co-content function~\cite{collins1978solving,dem:ofd}. The original formulations were later extended to include pressure-dependent demands and flow regulating devices ~\cite{deuerlein2009,moosavian2014hydraulic,Deuerlein2019}.  The advantages of optimization-based approaches are clear, linear and convex models can be efficiently solved to global optimality for very large networks using modern solvers \cite{mosek, gurobi}. The approach presented in this paper is most closely related to \cite{sela2015control} that initially proposed a geometric programming (GP) approximation ~\cite{duffin:gp} for solving the WFP by converting the nonconvex head loss equations into a GP form resulting in a nonlinear but convex optimization problem and, hence, a globally optimal solution is guaranteed. An important contribution of the previously proposed GP method is that it is non-iterative (i.e., a one-shot optimization problem). However, it is only applicable under the assumptions of a tree network topology, known and fixed flow directions, and was limited in the modeling complexities of valves and pumps. These assumptions make the previously proposed approach ~\cite{sela2015control} not suitable for urban water networks comprising branched and looped topologies.

In this study, we propose a novel GP approximation-based optimization approach to solve the network flow problem by taking advantage of the special form of the head loss equations. The main advantages and contributions of the proposed approach compared with previous GP-based modeling~\cite{sela2015control} are:  \textit{(1)} after transformation of variables and constraints the optimization problem that solves the WFP is linear, \textit{(2)} any arbitrary topologies and various types of valves and pumps can be seamlessly modeled providing modeling flexibility, and \textit{(3)} prior knowledge on flow directions or maximum flow rates is not required. 
Our approach involves two steps: \textit{(1)} the nonlinear nonconvex WFP is transformed to a nonlinear but convex problem using GP and \textit{(2)} the convex GP form is further transformed into linear form resulting in a set of linear equations.
In short, the proposed approach reduces the WFP to a system of linear equations and solves a series of linear programs (LP),  
thereby graciously scaling to large WDNs. Additionally, we demonstrate that the proposed approach can be straightforwardly extended to model pressure driven demands and leaks and integrated in control and optimization problems. We provide convergence proof, explore the sensitivy of the approach and propose acceleration scheme for computational speedup. The paper organization is given as follows. Section~\ref{sec:Model} describes the modeling of WDNs. Section~\ref{sec:gp-based} provides some necessary mathematical background related to geometric programming. Section~\ref{sec:GP} presents the paper's main contribution and the proposed algorithm and~\ref{sec:Convergence} presents the convergence proof. Section~\ref{sec:test} demonstrates the application of our approach to several case studies and Section~\ref{sec:sensitivity} presents the sensitivity analysis. Section~\ref{sec:future} proposes further extensions including pressure driven modeling and WFP-constrained optimization. Finally, Section~\ref{sec:conclusions} concludes the paper.

The following notations are used in the text -- italicized, boldface upper and lower case characters represent matrices and column vectors: $a$ is a scalar, $\m a$ is a vector, and $\m A$ is a matrix. Matrix $\m I$ denotes the identity square matrix, whereas $\m 0_{m \times n}$ denotes a zero matrix of with size $m$-by-$n$. {The notation $\mathbb{R}$ denotes the set of real numbers, and notations $\mathbb{R}^n$ and $\mathbb{R}^{m\times n}$ denote the sets of column vectors with $n$ elements and matrices with $m$-by-$n$ elements in $\mathbb{R}$.  For $\m x \in\mbb{R}^m$, $\m y \in \mbb{R}^n$, a compact column vector in $\mathbb{R}^{m+n}$  is defined as $\{\m x\, , \m y\} = [\m x^\top \ \m y^\top]^\top$. For matrices $\m A_{m_1 \times n}$ and $\m B_{m_2 \times n}$, the notation $[\m A; \m B]$ is defines as $[\m A^\top \ \m B^\top]^\top$. The element-wise product is represented as $\m x \circ \m y$ for $\m x, \m y \in\mbb{R}^m$.}
{The variables with upper case characters $\m \cdot^{\mathrm{J}}$, $\m \cdot^{\mathrm{R}}$, $\m \cdot^{\mathrm{TK}}$,  $\m \cdot^{\mathrm{P}}$, $\m \cdot^{\mathrm{M}}$, and  $\m \cdot^{\mathrm{W}}$ represent the variables related to junctions, reservoirs, tanks, pipes, pumps, and valves. }

\section{Modeling of WDNs}~\label{sec:Model}
A WDN is represented here by a directed graph {$\mathcal{G} = (\mathcal{V},\mathcal{E})$. The set $\mathcal{V}$ defines the nodes and is partitioned as $\mathcal{V} = \mathcal{J} \bigcup \mathcal{T} \bigcup \mathcal{R}$ where $\mathcal{J}$, $\mathcal{T}$, and $\mathcal{R}$ stand for the collection of $n_j$ junctions, $n_t$ tanks, and $n_r$  reservoirs, respectively. The set $\mathcal{E} \subseteq \mathcal{V} \times \mathcal{V}$ defines the links and is the partitioned as $\mathcal{E} = \mathcal{P} \bigcup \mathcal{M} \bigcup \mathcal{W}$, where $\mathcal{P}$, $\mathcal{M}$, and $\mathcal{W}$ represent the collection of $n_p$ pipes, $n_m$  pumps, and $n_w$ valves, respectively. The directed graph  $\mathcal{G}$ can be expressed by its  incidence matrix $\m A_\mathcal{G}$ which stands for the connection relationship between vertices and edges.} For the $i^\mathrm{th}$ node, the neighboring nodes are defined by the set $\mathcal{N}_i$, which is partitioned as $\mathcal{N}_i = \mathcal{N}_i^\mathrm{in} \bigcup \mathcal{N}_i^\mathrm{out}$, where $\mathcal{N}_i^\mathrm{in}$ and $\mathcal{N}_i^\mathrm{out}$  collect the nodes of the adjacent inflow and outflow links. Notice that the assignment of direction to each link (and the resulting inflow/outflow node classification) is arbitrary.  {Thus, $\m A_\mathcal{G}$ is comprised of $1$, $-1$, and $0$ elements indicating positive, negative, or no connection, respectively. $\m A_\mathcal{G}$ can be represented using the block column partition $[{\m A_{\m h}^\mathrm{P}}^\top {\m A_{\m h}^\mathrm{M}}^\top {\m A_{\m h}^\mathrm{W}}^\top]$, corresponding to pipe, pump, and valve edges, and block row partition ${[{\m A_{\m q}^\mathrm{J}}^\top {\m A_{\m q}^\mathrm{R}}^\top {\m A_{\m q}^\mathrm{TK}}^\top]^\top}$, corresponding to junction, reservoir, and tank nodes, as in~\eqref{equ:WFP-Incidence-matrix}. Note that the dimension of ${\m A_{\mathrm{P}}^\mathrm{J}}$ is $n_j \times n_p$, and the size of the other submatrices can be inferred similarly. The details of $\m A_\mathcal{G}$ are discussed in Section~\ref{sec:WFP}. Tab.~\ref{table:sets} summarizes the variables notation used in this paper.} 

\begin{figure*}
	{ 	\setlength\extrarowheight{5pt}
		\begin{align} ~\label{equ:WFP-Incidence-matrix}
		{\large \m A_\mathcal{G}} \hspace{-0.8ex}=\hspace{-0.8ex} \begin{array}{*{4}{cccc}@{}c}
		&\textit{Pipe}\ (n_p)&\textit{Pump} \ (n_m)&\textit{Valve} \ (n_w)\\
		\cline{2-4}
		\textit{Junction} \ (n_j)&\multicolumn{1}{|c|}{{\m A_{\mathrm{P}}^\mathrm{J}}}&\multicolumn{1}{|c|}{{\m A_{\mathrm{M}}^\mathrm{J}}}&\multicolumn{1}{|c|}{{\m A_{\mathrm{W}}^\mathrm{J}}}& \multirow{1}{*}{$\left.\rule[1.ex]{0pt}{1.ex}\right\}\m A_{\m q}^\mathrm{J}$}\\
		\cline{2-4}
		\textit{Reservoir} \ (n_r)&\multicolumn{1}{|c|}{\m A_{\mathrm{P}}^\mathrm{R}}&\multicolumn{1}{|c|}{\m A_{\mathrm{M}}^\mathrm{R}}&\multicolumn{1}{|c|}{\m A_{\mathrm{W}}^\mathrm{R}}&  \multirow{1}{*}{$\left.\rule[1.ex]{0pt}{1.ex}\right\}\m A_{\m q}^\mathrm{R}$}\\ 
		\cline{2-4}
		\textit{Tank} \ (n_t)&\multicolumn{1}{|c|}{\m A_{\mathrm{P}}^\mathrm{TK}}&\multicolumn{1}{|c|}{\m A_{\mathrm{M}}^\mathrm{TK}}&\multicolumn{1}{|c|}{\m A_{\mathrm{W}}^\mathrm{TK}}&\multirow{1}{*}{$\left.\rule[1.ex]{0pt}{1.ex}\right\}\m A_{\m q}^\mathrm{TK}$}\\ 
		\cline{2-4}
		\noalign{\vspace{-6pt}}
		\multicolumn{1}{c}{}	& \multicolumn{1}{@{}c@{}}{\underbrace{\hspace*{1\tabcolsep}\hphantom{.......}}_{{\Large {\m A_{\m h}^\mathrm{P}}^\top}}}  & \multicolumn{1}{@{}c@{}}{\underbrace{\hspace*{1\tabcolsep}\hphantom{.......}}_{{\Large {\m A_{\m h}^\mathrm{M}}^\top}}} & \multicolumn{1}{@{}c@{}}{\underbrace{\hspace*{1\tabcolsep}\hphantom{.......}}_{{\Large {\m A_{\m h}^\mathrm{W}}^\top}}}
		\end{array}
		\end{align}}
	\hrulefill
\end{figure*}
\normalcolor
\begin{table*}[t]
	\caption{Variable notation.}
	\small
	\tabcolsep= 0.7cm
	\renewcommand{\arraystretch}{1.2}
	\centering
	\begin{tabular}{ cc }
		\hline
		\textit{Notation} & \textit{Description} \\ \hline
		$h_i^{\mathrm{J}}$, $h_i^{\mathrm{R}}$, $h_i^{\mathrm{TK}}$ &  Head at the $i^\mathrm{th}$ junction, reservoir, or tank   \\ 
		{$q_{ij}^{\mathrm{P}}$}, {$q_{ij}^{\mathrm{M}}$}, {$q_{ij}^{\mathrm{W}}$} & {Flow through the pipe, pump, or valve from node $i$ to node $j$}  \\
		$\Delta h_{ij}^{\mathrm{P}}$, $\Delta h_{ij}^{\mathrm{M}}$, $\Delta h_{ij}^{\mathrm{W}}$&  Head loss or gain  from $i$ to $j$ for the pipe, pump, or valve \\ 
		{$s_{ij}$} & {Speed of the pump through node $i$ to node $j$}  \\	
		{$o_{ij}$} & {Openness of the valve through node $i$ to node $j$}  \\
		$\m {\xi}$ & {A vector collecting all variables (head and flow)} \\
		$ \hat{\m \xi}$ & {GP form of $\m \xi$}\\ 
		{$\langle{\hat{\m \xi}}\rangle_{n}$} & {The $n^\mathrm{th}$ iteration value of $\hat{\m \xi}$}  \\	
		$\m \xi_{\mathrm{EPANET}}$ & Solution provided by EPANET software. \\ 
		$ \m \xi_{\mathrm{GP-LP}}$ & Solution from our proposed GP-LP-based approach. \\ 
		\hline		
	\end{tabular}
	\label{table:sets}
\end{table*}
\subsection{Modeling components}~\label{sec:Model_iass}
The basic hydraulic equations describing the flow in WDNs are derived from the principles of {conservation of mass} and {energy} ~\cite{todini2013}. {For elements such as nodes, conservation of mass means the sum of inflows and outflows is equal to zero, and for storage tanks to the change in the water storage volume.} The conservation of energy states that the energy difference stored in a component is equal to the energy increases minus energy losses, such as, frictional and minor losses~\cite{puig2017real}. According to these basic laws, the equations that model mass and energy conservation for all components in WDNs can be written in explicit and compact matrix-vector forms, as detailed next.
\subsubsection{Tanks and reservoirs}	
We assume that reservoirs have infinite water supply and the head of the $i^\mathrm{th}$ reservoir is fixed~\cite{zamzam2018optimal,singh2018optimal,gleixner2012towards} 
and we have
\begin{linenomath*}
	\begin{equation}\label{equ:reservoir}
	h_i^{\mathrm{R}} = h_i^{\mathrm{R}_\mathrm{set}},
	\end{equation} 
\end{linenomath*}
where $h_i^{\mathrm{R}_\mathrm{set}}$ is specified. 

The head created by a cylindrical tank that has a fixed cross sectional area can be described as
\begin{linenomath*}
	\begin{equation}\label{equ:tank}
	h_i^{\mathrm{TK}} = h_i^{\mathrm{TK}_\mathrm{set}},
	\end{equation} 
\end{linenomath*}
where $h_i^{\mathrm{TK}_\mathrm{set}} = \frac{V_i}{A_i^{\mathrm{TK}}} + E_i^{\mathrm{TK}}$ and the elevation $E_i^{\mathrm{TK}}$,  volume $V_i$ and cross sectional area $A_i^{\mathrm{TK}}$ of the $i^\mathrm{th}$ tank can be measured.

\subsubsection{Junctions and pipes}
Junctions are points of connection between links where water flow merges or splits. The expression of mass conservation of the $i^\mathrm{th}$ junction  can be written as
\begin{linenomath*}
	\begin{equation}\label{equ:nodes}
	\sum_{j \in \mathcal{N}_i^\mathrm{in}} q_{ji} - \sum_{j \in \mathcal{N}_i^\mathrm{out}} q_{ij} = d_i,
	\end{equation} 
\end{linenomath*}
where $d_i$ stands for end-user demand that is extracted from node $i$, and we assume that the demand is known for the WFP. The major head loss of a pipe from node $i$ to $j$ is due to friction and is determined by 
\begin{linenomath*}
	\begin{equation}~\label{equ:head-flow-pipe}
	\Delta h_{ij}^\mathrm{P}  = h_{i} - h_{j} = R_{ij} {q_{ij}^\mathrm{P}}|q_{ij}^\mathrm{P}|^{\mu-1},
	\end{equation} 
\end{linenomath*}
{where $R_{ij}$ is pipe resistance coefficient, which is a function of pipe size, length, and material; $\mu$ is the constant flow exponent. Note that $R_{ij}$ and $\mu$ vary correspondingly with the most common formulae to model the head loss, which are  \textit{Hazen-Williams}, \textit{Darcy-Weisbach}, and \textit{Chezy-Manning}~\cite{linsley1979water,rossman2000epanet}. The approach presented in this paper considers any of the three formulae.  The minor head losses in pipes caused by turbulence that occurs at bends and fittings are not considered in this paper, but could be easily modeled using surrogate pipe length.}

\subsubsection{Pumps} 
A head increase/gain can be generated by a pump between the suction node $i$ and the delivery node $j$. The pump properties dictate the relationship function between the pump flow and head increase~\cite{linsley1979water}. Generally, the head gain can be expressed as
\begin{linenomath*}
	\begin{align} \label{equ:head-flow-pump}
	\Delta h_{ij}^\mathrm{\mathrm{M}} = h_{i} - h_{j} = -{s_{ij}^2}\left(h_0 - r  (q_{ij}^\mathrm{M} s_{ij}^{-1})^\nu \right),
	\end{align} 
\end{linenomath*}
where $h_0$ is the shutoff head, $q_{ij}^\mathrm{M} $ is the flow, $s_{ij} \in (0,s_{ij}^{\mathrm{max}} ]$ is the relative speed, {which is known, $r$ and $\nu$ are the curve coefficients of the pump that are chosen from a particular range of values}. It is worthwhile to notice that the head gain $h_{ij}^M$ is always a negative value and the flow through the pump is always strictly positive. Pump flow and head constraints will be later modeled as operational constraints~\eqref{equ:headgainlimit} and~\eqref{equ:flowLimit}. 

\subsubsection{Valves} 
Several types of valves can be utilized to regulate the flows or pressures in WDNs.
General Purpose Valves (GPV), Pressure Reducing Valves (PRV), and Flow Control Valves (FCV) are commonly used valves that are controlled through valve openness or set points for pressure reduction or flow regulation. The different valve flow-head relationships used in our paper are based on ~\cite[Chapter 3]{rossman2000epanet}.
{GPVs} can be used to model turbines, well draw-down or reduced-flow backflow prevention valves. 
Here, we assume that the GPVs are modeled similarly to a pipe with controlled resistance coefficient, which  can be expressed as
\begin{linenomath*}
	\begin{align}~\label{equ:head-flow-valve}
	\Delta  h_{ij}^\mathrm{W}  = h_{i} - h_{j} = o_{ij}^{-1} R_{ij} {q_{ij}^\mathrm{W}}|q_{ij}^\mathrm{W}|^{\mu-1}, 
	\end{align} 
\end{linenomath*}
where {$o_{ij} \in (0,1]$} is {a known parameter} depicting the openness of a valve, and the rest of the variables are similar to the pipe model. When $o_{ij} = 1$ the valve is fully open and as $o_{ij}$ decreases, the valve closes resulting in greater losses~\cite{piller2014modeling}. 
When a GPV is completely closed, no constraint exists between $h_i $ and $h_j$ indicating that the two nodes are decoupled and the corresponding constraint~\eqref{equ:head-flow-valve} should be removed.  

{PRVs} limit the pressure at a specific location in the network (reverse flow is not allowed) and set the pressure to $P^\mathrm{set}$ on its downstream side when the upstream pressure is higher than $P^\mathrm{set}$~\cite[Chapter 3.1]{rossman2000epanet}, {otherwise, they are treated as open pipes with minor head loss}. Assuming that the upstream side is denoted as $i$, and the downstream side is $j$ and given the status of a PRV, the PRV can be modeled as 
\begin{linenomath*}
	\begin{subnumcases}{~\label{equ:head-prv-valve}}
	\Delta h_{ij}^\mathrm{W} = h_{i} - h_{j} = l_{ij} {q_{ij}^\mathrm{W}}|q_{ij}^\mathrm{W}|,\mathrm{OPEN}\\
	h_{j}  = h^{\mathrm{W}_\mathrm{set}}  ,\;\text{ACTIVE},~\label{equ:head-prv-valve-active}
	\end{subnumcases}
\end{linenomath*}
where {$l_{ij}$ is the lumped minor head loss coefficient depending on the acceleration of gravity, cross-sectional area, and local losses of the PRV}. Parameter $h^{\mathrm{W}_\mathrm{set}}  $ is the pressure setting converted to head implying $h^{\mathrm{W}_\mathrm{set}}   = E_j + P^\mathrm{set}$, and $E_j$ is the elevation at junction $j$, parameter $P^\mathrm{set}$ is the pressure setting of the PRV and both are constant. Therefore, the head $h_j$ is fixed, and the fact that reverse flow is not allowed in PRVs can be expressed as operational constraint included in~\eqref{equ:flowLimit}. Similarly for GPV, if PRV completely closed, the constraint between the two adjacent nodes $h_i $ and $h_j$ is removed. 

{{FCVs}  limit the flow to a specified setting  $q^{\mathrm{W}_\mathrm{set}}  $ when the head $h_{i}$ at upstream node~$i$ is greater than the head $h_{j}$ at downstream node~$j$, otherwise, FCVs are treated as open pipes with minor head loss. FCVs can be modeled as}
\begin{linenomath*}
	\begin{subnumcases}{~\label{equ:head-fcv-valve}}
	\Delta h_{ij}^\mathrm{W} = h_{i} - h_{j} = l_{ij} {q_{ij}^\mathrm{W}}|q_{ij}^\mathrm{W}|,\mathrm{OPEN}\\
	q^{\mathrm{W}}_{ij} = q^{\mathrm{W}_\mathrm{set}}  ,\;\text{ACTIVE}, ~\label{equ:head-fcv-valve-active}
	\end{subnumcases}
\end{linenomath*}
where  {$l_{ij}$ is the lumped minor head loss coefficient and $q^{\mathrm{W}_\mathrm{set}}  $ is the setting value}. 

{Note that regardless of the type of link, i.e. pipe, pump, or valve, a closed indicates that the corresponding flow $q_{ij}$ is zero, thus the corresponding links are removed from the incidence matrix $\m A_\mathcal{G}$, and inherently no constraints are imposed between its adjacent nodes $i$ and $j$.}

\subsection{Nonlinear water flow problem formulation}\label{sec:WFP}
This section derives an optimization-based formulation given the WDN model. Aside from the physical constraints listed above, typical design and operation problems pertaining to WDNs also consider engineering constraints, such as restricting the desired flows and heads in the network. These additional constraints can be written as
\begin{linenomath*}
	\begin{subequations} ~\label{equ:constraints} 
		\begin{align}
		h_{i}^{\mathrm{min}} &\leq  h_{i} \leq h_{i}^{\mathrm{max}}~\label{equ:tankLimit} \\
		q_{ij}^{\mathrm{min}} &\leq  q_{ij} \leq q_{ij}^{\mathrm{max}}~\label{equ:flowLimit}\\
		&h_{ij}^{\mathrm{M}} \leq 0.~\label{equ:headgainlimit}
		\end{align}
	\end{subequations}
\end{linenomath*}
Eqs.~\eqref{equ:tankLimit}--\eqref{equ:flowLimit} are the lower and upper bounds on the heads of nodes, flows through links; Eq.~\eqref{equ:headgainlimit} is the head increase delivered by pumps. Let the compact vectors $\m h^{\mathrm{J}}$, $\m h^{\mathrm{R}}$, and $\m h^{\mathrm{TK}}$ collect the heads at junctions, reservoirs, and tanks, { $\m h  \triangleq  \lbrace \m h^{\mathrm{J}}, \m h^{\mathrm{R}}, \m h^{\mathrm{TK}}   \rbrace$ collect all the heads at the nodes, where $\m h \in \mathbb{R}^{n_h}$ and $n_h = n_j + n_r + n_t$ is the summation of} the number of junction, reservoirs, and tanks, respectively. Similarly, the flow through pipes, pumps, and valves are collected by compacted vectors $\m q^{\mathrm{P}}$, $\m q^{\mathrm{M}}$, and $\m q^{\mathrm{W}}$, {let  $\m q  \triangleq  \lbrace \m q^{\mathrm{P}}, \m q^{\mathrm{M}}, \m q^{\mathrm{W}}   \rbrace$, and $\m q \in \mathbb{R}^{n_q}$, where $n_q = n_p + n_m + n_w$ is the summation of} the number of pipes, pumps, and valves, respectively. We define a vector collecting all above optimization variables  as $\m \xi  \triangleq  \lbrace \m h, \m q  \rbrace$, and  $\m \xi \in \mathbb{R}^{n_{\xi}}$ where $n_{\xi} = n_h + n_q$. Thus, all constraints can be summarized as $\m \xi \in [\m \xi^{\mathrm{min}},\m \xi^{\mathrm{max}}]$, and Eqs.~\eqref{equ:reservoir}-\eqref{equ:head-fcv-valve} can be presented as
\begin{linenomath*}
	\begin{equation}\label{equ:nonlinearMatrix}
	\hspace{-1em}\includegraphics[width=0.920\linewidth,valign=c]{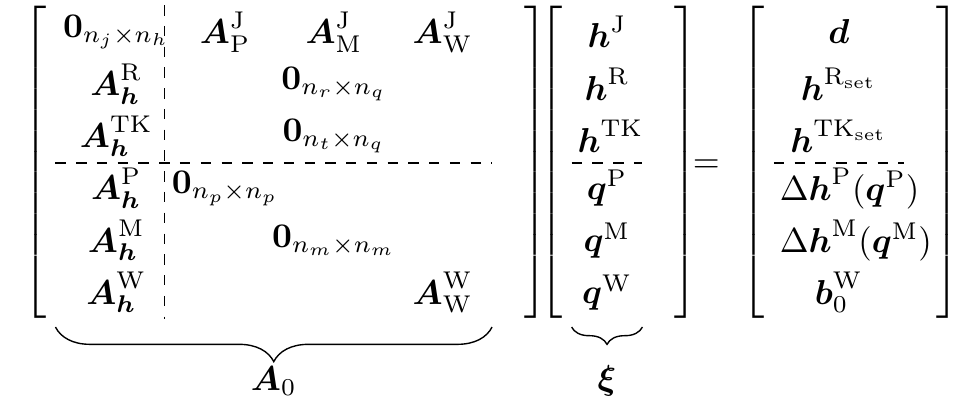},
	\end{equation}
\end{linenomath*}
where ${\m A}_{\m q}^{\mathrm{J}} \triangleq [ \m A_{\mathrm{P}}^\mathrm{J}\ \m A_{\mathrm{M}}^\mathrm{J}\  \m A_{\mathrm{W}}^\mathrm{J} ]$ is the mass balance for all nodes~\eqref{equ:nodes}, ${\m A}_{\m h}^{\mathrm{R}} $ and ${\m A}_{\m h}^{\mathrm{TK}}$ collect~\eqref{equ:reservoir} and~\eqref{equ:tank} for reservoirs and tanks,  ${\m A}_{\m h}^{\mathrm{P}}$ and ${\m A}_{\m h}^{\mathrm{M}}$ collect the head loss equation~\eqref{equ:head-flow-pipe} and the head gain equation~\eqref{equ:head-flow-pump}, and ${\m A}_{\m h}^{\mathrm{W}}$ and ${\m A}_{\mathrm{W}}^{\mathrm{W}}$ are  the left hand side of valve~\eqref{equ:head-flow-valve},~\eqref{equ:head-prv-valve}, and~\eqref{equ:head-fcv-valve}.
The right hand side  is a  vector collecting the corresponding coefficients including demand $\m d \in \mathbb{R}^{n_{j}}$ for all junctions,  settings of reservoirs ${\m h}^{\mathrm{R}_\mathrm{set}} \in \mathbb{R}^{n_r}$ and tanks ${\m h}^{\mathrm{TK}_\mathrm{set}} \in \mathbb{R}^{n_t}$, {nonlinear head loss/gain functions for all links ${\Delta \m h}^{\mathrm{P}}({\m q}^{\mathrm{P}})$, ${\Delta \m h}^{\mathrm{M}}({\m q}^{\mathrm{M}})$, and ${\Delta \m h}^{\mathrm{W}}({\m q}^{\mathrm{W}})$ or valve settings $\m h^{\mathrm{W}_\mathrm{set}}$ and $\m q^{\mathrm{W}_\mathrm{set}}$ collected in $\m b_0^\mathrm{W}$.}

\begin{remark}
	$\m A_0$ is a square matrix, and  $\m A_0 \in \mathbb{R}^{n_{\xi} \times n_{\xi}} $ since  ${\m A}_{\m q}^{\mathrm{J}} \in \mathbb{R}^{n_{j} \times n_q}$, ${\m A}_{\m h}^{\mathrm{R}} \in \mathbb{R}^{n_{r} \times n_h}$, ${\m A}_{\m h}^{\mathrm{TK}} \in \mathbb{R}^{n_{t} \times n_h}$, ${\m A}_{\m h}^{\mathrm{P}} \in \mathbb{R}^{n_{p} \times n_h}$, ${\m A}_{\m h}^{\mathrm{M}} \in \mathbb{R}^{n_{m} \times n_h}$, and ${\m A}_{\m h}^{\mathrm{W}} \in \mathbb{R}^{n_{w} \times n_h}$.
\end{remark}

\begin{assumption}\label{asp:1}
	$\m A_0$ is invertible.
\end{assumption}
\normalcolor
The overall nonlinear modeling of \textsc{WFP} can be written as
\begin{linenomath*}
	\begin{align}
	\textsc{\textbf{WFP:}}\;\;\;\;\;   \mathrm{find} \;\; \;\;& \hspace{2em}{\boldsymbol \xi}  \label{equ:WFP} \\
	\mathrm{s.t.}\;\;\;\;& ~\eqref{equ:constraints}, \eqref{equ:nonlinearMatrix} .\notag
	\end{align}
\end{linenomath*}
The WFP~\eqref{equ:WFP} is nonlinear and nonconvex due to the head loss models of pipes and pumps. 

Motivated by the literature gaps discussed in Section~\ref{sec:literature}~\cite{sela2015control}, we propose a new GP-based optimization approach to solve the WFP, which is convex in the variables, considers various kinds of valves and pumps, while not requiring the a priori knowledge of water flow direction, and applies to any network topology.  Note that the GP form itself is not convex, however, the $\log$ form of GP is convex. Hence, when we say GP is convex in this paper, we mean that the $\log$ form of GP is convex.  After the GP is developed, the problem is transformed to an LP, as discussed in Section~\ref{sec:LP-GP}.

\section{Geometric programming and modeling tricks}~\label{sec:gp-based}
A basic introduction to GP is given in this section and a simple LP example is presented to illustrate how to convert a problem into its GP form.
\subsection{Geometric programming}~\label{sec:GPmodeling}
A geometric program is a type of optimization problem with objective and constraint functions that are monomials and posynomials~\cite{boyd2007tutorial}.
{A real valued function $g(\m x) = c x_1^{a_1} x_2^{a_1} \cdots x_n^{a_n}$, where $c > 0$, $\m x>0$, and $a_i \in \mathbb{R}$, is called a \textit{monomial} of the variables $x_1, \cdots, x_n$}. A sum of one or more monomials, i.e., a function of the form $f(\m x) = \sum_{k=1}^{K}c_k x_1^{a_{1k}} x_2^{a_{2k}} \cdots x_n^{a_{nk}}$
where $c_k > 0$, is called a \textit{posynomial} with $K$ terms in the vector variable $\m x$.
A standard GP can be written as
\begin{linenomath*}
	\begin{align}~\label{equ:GP-standard}
	\textit{GP:}\;\; \min_{\m x >0 } \hspace{15pt} &f_0(\m x) \notag  \\
	\mathrm{s.t.}\hspace{15pt}& f_i(\m x) \leq 1, i = 1,\cdots, m \\
	& g_i(\m x) = 1, i = 1,\cdots, p, \notag 
	\end{align} 
\end{linenomath*}
where $\m x$ is an entry-wise positive optimization variable, $f_i(\m x)$ are posynomial functions and $g_i(\m x)$ are monomials. A standard GP form is nonlinear and nonconvex. The main technique to solving a GP efficiently is to convert it to a nonlinear but convex problem using a logarithmic change of variables, and a logarithmic transformation of the objective and constraint functions. Defining $y_i = \log x_i$, the $\log$ form of GP can be expressed as
\begin{linenomath*}
	\begin{align*}
	\textit{GP-$\log$:}\;\; \min_{y } \hspace{15pt} &\log f_0(\m y) \notag  \\
	\mathrm{s.t.}\hspace{15pt}& \log f_i(\m y) \leq 0, i = 1,\cdots, m \\
	& \log g_i(\m y) = 0, i = 1,\cdots, p, \notag 
	\end{align*} 
\end{linenomath*}
The $\log$ form of GP is convex, and can be solved efficiently with modern solvers~\cite[Section 2.4]{boyd2007tutorial}. We note that the  $\log$ function can be  with any base $b$ which is greater than $1$.  With above analysis, we can see that our task is to formulate our WFP as a standard GP, and then the rest can be solved directly by GP solver ~\cite{mosek,cvx}. 

\subsection{{Handling negative flows}}
Since the direction of the flows in the WFP problem is generally unknown,  
the flow in each pipe can be viewed as free variable, i.e., it is unrestricted in its sign. However, all variables in \eqref{equ:GP-standard} are required to be positive, and reverse direction of flows can not be directly modeled using GP. Several techniques are available to convert the free variables to positive variables, for example, by representing the positive and negative parts by positive dummy variables as in the simplex method~\cite{lustig1994interior} for solving LPs. However, this modeling trick does not apply in our setting (see ~\eqref{equ:head-flow-pipe} and ~\eqref{equ:head-flow-valve}). Here, we propose another trick to convert a free variable to a positive one, thus allowing modeling reverse direction of flows. Consider  an exponential function $f(x) = b^x$ that can map any $x$ to $f(x)$ and $f(x)$ is always positive. Taking advantage of this idea, we can convert a problem with a negative feasible region into a new problem with a positive feasible region. After the solution of the transformed problem is obtained, the original solution can be obtained simply by reverting back. For ease of demonstration, we illustrate this idea using a simple LP problem~\eqref{equ:LP}.

First, the original free variables are converted into the corresponding GP variables denoted using $\hat{x}$, e.g., the variables $x_1$ and $x_2$ turn into $\hat{x}_1$ and $\hat{x}_2$ via $\hat{x}_1 = b^{x_1}$ and $\hat{x}_2 = b^{x_2}$ where the base $b > 1$.  Second, the constraints and objective functions are converted into the monomial or posynomial form, e.g., {the constraint $ -x_1 + x_2 \leq 4$ can be expressed as $b^{-4}\;\hat{x}_1^{-1}\; \hat{x}_2^1  \leq 1$ via executing exponential on both sides of $ -x_1 + x_2 \leq 4$}. Thus, we successfully convert an LP~\eqref{equ:LP} into its GP form~\eqref{equ:GP-linear} and the solution of the original LP problem can also be obtained using GP solver.
\begin{multicols}{2}
	\noindent 	\begin{linenomath*}\begin{align}
		\min \ \ & 2x_1 + 3x_2 ~\label{equ:LP} \\
		\mathrm{s.t.}\;& -x_1 + x_2 \leq 4 \notag
		\end{align} 
	\end{linenomath*}
	\begin{align}
	\min \;& \hat{x}_1^2\; \hat{x}_2^3   ~\label{equ:GP-linear}\\
	\mathrm{s.t.}\;&   b^{-4}\;\hat{x}_1^{-1}\; \hat{x}_2^1  \leq 1 \notag
	\end{align}
\end{multicols}
Any LP problem can be converted, however, the conversion of LP seems to make the transformed problem harder. This is not the case when the nominal problem is highly nonlinear and nonconvex. The technique we introduced may help to transform a nonconvex, nonlinear problem into a convex one as shown in subsequent sections.

\section{GP modeling and corresponding LP modeling of WFP}~\label{sec:GP}
Based on the new introduced optimization technique, we convert the nonlinear WFP ~\eqref{equ:WFP} into its GP form, then derive the corresponding LP model, and propose an algorithm to solve the GP-LP problem.
\subsection{Conversion of variables} 
Here, the GP variables $\hat{\boldsymbol \xi}$ are obtained by mapping the optimization variables $\m \xi$ in~\eqref{equ:WFP}. 
Specifically, we convert the head and demand at the $i^\mathrm{th}$ node, $h_i$ and $d_i$, and the flow $q_{ij}$ into positive values ${\hat{h}_i}$, ${\hat{d}_i}$, and ${\hat{q}_{ij}}$ through exponential functions, as follows
\begin{linenomath*}
	\begin{equation}~\label{equ:NLPGP}
	{\hat{h}_i} \triangleq {b}^{h_i}, \; {\hat{d}_i} \triangleq {b}^{d_i},\;
	{\hat{q}_{ij}} \triangleq {b}^{q_{ij}} , \; 
	\end{equation} 
\end{linenomath*}
where $b=1+\delta$ is a constant base and $\delta$ is a small positive number. The variables ${\hat{h}_i}$, ${\hat{d}_i}$, and ${\hat{q}_{ij}}$ are positive which can then be used to transform the nonconvex WFP~\eqref{equ:WFP} into a GP. 
\subsection{ Conversion of mass and energy balance equations}~\label{sec:modelconversion}
\subsubsection{{Mass balance equations for junctions}}\label{sec:massbalance}
{Converting mass balance at junctions following the above exponential mapping~\eqref{equ:NLPGP} is straightforward. After exponentiating both sides of~\eqref{equ:nodes}, variables collected in $\m \xi$ are changed into ${\hat{\m \xi}}$, the summation is turned into multiplication, and we obtain
	\begin{linenomath*}
		\begin{align*}
		b^{\sum_{j \in \mathcal{N}_i^\mathrm{in}} q_{ji} - \sum_{j \in \mathcal{N}_i^\mathrm{out}} q_{ij}} &= \prod_{j \in \mathcal{N}_i^\mathrm{in}}\hspace{-3pt}b^{{q}_{ji}}\hspace{-6pt} \prod_{j \in \mathcal{N}_i^\mathrm{out}}\hspace{-3pt}{b^{-{q}_{ij}}} \\
		= \prod_{j \in \mathcal{N}_i^\mathrm{in}}\hspace{-3pt}\hat{q}_{ji}\hspace{-6pt} \prod_{j \in \mathcal{N}_i^\mathrm{out}}\hspace{-3pt}{\hat{q}_{ij}}^{-1}&= b^{d_i}=\hat{d}_i.
		\end{align*} 
	\end{linenomath*}
	After the transformation, constraint~\eqref{equ:nodes} is converted to monomial equality constraint written as
	\begin{linenomath*}
		\begin{align}
		\prod_{j \in \mathcal{N}_i^\mathrm{in}}\hspace{-3pt}\hat{q}_{ji}\hspace{-6pt} \prod_{j \in \mathcal{N}_i^\mathrm{out}}\hspace{-3pt}{\hat{q}_{ij}}^{-1}{{\hat{d}_i}^{-1}} &= 1.  ~\label{equ:nodes-exp}
		\end{align} 
	\end{linenomath*}
	
	\normalcolor
	
	\subsubsection{{Energy balance equations for pipes}}
	Now we convert the head loss model for pipes, and let ${\Delta \hat{h}_{ij}^{\mathrm{P}}}$ be the GP form of head loss of a pipe, which is obtained  by exponentiating both sides of~\eqref{equ:head-flow-pipe} as follows
	\begin{linenomath*}
		\begin{align*}
		{\hat{h}_{i}} {\hat{h}_{j}^{-1}} &= {\Delta  \hat{h}_{ij}^{\mathrm{P}}}  = {b^{\left(q_{ij}^{\mathrm{P}} R_{ij} {|q_{ij}^{\mathrm{P}}|}^{\mu-1} - q_{ij}^{\mathrm{P}} + q_{ij}^{\mathrm{P}}\right)}}\\
		&= b^{q_{ij}^{\mathrm{P}} \left(R_{ij} {|q_{ij}^{\mathrm{P}}|}^{\mu-1} - 1\right)}\ {\hat{q}_{ij}}= \hat{c}^{\mathrm{P}}(q_{ij}^{\mathrm{P}})\ {\hat{q}_{ij}^{\mathrm{P}}},
		\end{align*}
	\end{linenomath*}
	where  $\hat{c}^{\mathrm{P}}(q_{ij}^{\mathrm{P}})=b^{q_{ij}^{\mathrm{P}} \left(R {|q_{ij}^{\mathrm{P}}|}^{\mu-1} - 1\right)}$ is a function of $q_{ij}^{\mathrm{P}}$ {which means $\hat{c}^{\mathrm{P}}(q_{ij}^{\mathrm{P}})$ can be viewed as a known when $q_{ij}^{\mathrm{P}}$ is given.}
	{At first, we can make an initial guess denoted by $\langle{q_{ij}^{\mathrm{P}}}\rangle_0$ for the $0^\mathrm{th}$ iteration ($\langle{\hat{c}^{\mathrm{P}}}\rangle_0$ can be obtained if $\langle{q_{ij}^{\mathrm{P}}}\rangle_0$ is known), thus, for the $n^\mathrm{th}$ iteration, the corresponding values are denoted by $\langle{q_{ij}^{\mathrm{P}}}\rangle_n$ and $\langle{\hat{c}^{\mathrm{P}}}\rangle_n$. If the flow rates are close to each other between two successive iterations, we can approximate $\langle{\hat{c}^{\mathrm{P}}}\rangle_n$ using  $\langle{\hat{c}^{\mathrm{P}}}\rangle_{n-1}$, that is
		$\langle{\hat{c}^{\mathrm{P}}}\rangle_n \approx \langle{\hat{c}^{\mathrm{P}}}\rangle_{n-1}.$
		Then, for each iteration $n$,}  
	\begin{linenomath*}
		\begin{equation*}~\label{equ:hat-cp-exp}
		\langle{\hat{c}^{\mathrm{P}}}\rangle_n = b^{\langle{q_{ij}^{\mathrm{P}}}\rangle_{n-1} \left(R {|\langle{q_{ij}^{\mathrm{P}}}\rangle_{n-1}|}^{\mu-1} - 1\right)}
		\end{equation*} 
	\end{linenomath*}
	can be approximated given the flow value $\langle{q_{ij}^{\mathrm{P}}}\rangle_{n-1}$ from the previous iteration. With this approximation, the head loss constraint for each pipe can be written as a monomial equality constraint
	\begin{linenomath*}
		\begin{equation}~\label{equ:head-loss-pipe-exp}
		{\hat{h}_{i}} {\hat{h}_{j}^{-1}} [\hat{c}^{\mathrm{P}}]^{-1} [\hat{q}_{ij}^{\mathrm{P}}]^{-1} = 1.
		\end{equation} 
	\end{linenomath*}
	
	{The idea is to iteratively update the above monomial equality constraint, where the highly nonlinear term is included into a parameter $\hat{c}^{\mathrm{P}}$ and computed based on the solution of the previous iteration. The new obtained solution is used to update $\hat{c}^{\mathrm{P}}$ again and  generate the constraints in next iteration. This technique is similar to the iterative update in the gradient and Newton-Raphson approaches ~\cite{todini1987s}.}
	\subsubsection{{Energy balance equations for pumps}}
	Similarly, the new variables ${\hat{q}_{ij}^{\mathrm{M}}} = b^{q_{ij}^{\mathrm{M}}}$ and ${\hat{s}_{ij}} = b^{s_{ij}}$ for $(i,j) \in \mathcal{M}$ are introduced for pumps.  Let ${\Delta \hat{h}_{ij}^\mathrm{M}} $ be the GP form of head increase of a pump:
	\begin{linenomath*}
		\begin{align} \label{equ:pumpheadloss}
		{\hat{h}_{i}} {\hat{h}_{j}^{-1}} &= {\Delta \hat{h}_{ij}^\mathrm{M}} = b^{-{s_{ij}^2}(h_0 - r\; (q_{ij}^{\mathrm{M}})^{\nu}s_{ij}^{-\nu})} \\
		&=b^{-{s^2_{ij}} h_0}\ (b^{q_{ij}^{\mathrm{M}}})^{ r (q_{ij}^{\mathrm{M}})^{\nu-1} s_{ij}^{2-\nu}} ={\hat{c}_1^{\mathrm{M}}} ({\hat{q}_{ij}})^{c_2^{\mathrm{M}}}, \notag
		\end{align} 
	\end{linenomath*}
	where $\hat{c}_1^{\mathrm{M}} = b^{-{s_{ij}} h_0}$ and $c_2^{\mathrm{M}} = r (q_{ij}^{\mathrm{M}})^{\nu-1} s_{ij}^{2-\nu}$. Parameters $\hat{c}_1^{\mathrm{M}}$ and $c_2^{\mathrm{M}}$ follow a similar iterative process as $\hat{c}^{\mathrm{P}}$. That is, they are treated at the $n^\mathrm{th}$ iteration as constants based on the flow and relative speed values at the ${n\hspace{-2pt}-\hspace{-2pt}1}^\mathrm{th}$ iteration. Hence, the approximating equation for the pump head increase becomes the monomial equality constraint				
	\begin{linenomath*}
		\begin{equation}  \label{equ:head-flow-pump-exp}
		{\hat{h}_{i}} {\hat{h}_{j}^{-1}}[\hat{c}_1^{\mathrm{M}}]^{-1}[{\hat{q}_{ij}^{\mathrm{M}}}]^{-c_2^{\mathrm{M}}}  = 1.
		\end{equation} 
	\end{linenomath*}
	\subsubsection{{Energy balance equations for valves}}
	As for valves, the derivation of GPVs is the same as for pipes except for an extra variable ${\hat{o}_{ij}} = b^{o_{ij}^{-1}}$ for $(i,j) \in \mathcal{W}$ is introduced.  Let $\Delta {\hat{h}_{ij}^{\mathrm{W}}}$ be the GP form of head loss of a valve, which is obtained  by exponentiating both sides of~\eqref{equ:head-flow-valve} as follows.
	\begin{linenomath*}
		\begin{align*} 
		{\hat{h}_{i}} {\hat{h}_{j}^{-1}} &= {\hat{h}_{ij}^{\mathrm{W}}}  = { b^{\left(o_{ij}^{-1}q_{ij}^\mathrm{W} R {|q_{ij}^\mathrm{W}|}^{\mu-1} - q_{ij}^\mathrm{W} + q_{ij}^\mathrm{W}\right)}}\\
		&= b^{o_{ij}^{-1} \left(R q_{ij}^\mathrm{W} {|q_{ij}^\mathrm{W}|}^{\mu-1} - q_{ij}^\mathrm{W}\right)} \ {\hat{q}_{ij}} = {\hat{c}^{\mathrm{W}}}\ {\hat{q}_{ij}}, \notag
		\end{align*} 
	\end{linenomath*}
	where  ${\hat{c}^{\mathrm{W}}(q_{ij}^\mathrm{W})=b^{o_{ij}^{-1}q_{ij}^\mathrm{W} \left( R{|q_{ij}^\mathrm{W}|}^{\mu-1} - 1\right) }}$ is a similar parameter as the parameters in pipe and pump models. Hence,  the monomial equality constraint can be used for GPVs 			
	\begin{linenomath*}
		\begin{equation} ~\label{equ:head-loss-valve-gpv-exp}
		{\hat{h}_{i}} {\hat{h}_{j}^{-1}} [\hat{c}^{\mathrm{W}}]^{-1} [\hat{q}_{ij}^{\mathrm{W}}]^{-1} = 1.
		\end{equation} 
	\end{linenomath*}
	
	For PRVs and FCVs, the conversion process is similar as the one of pipes or GPVs, and Eqs.~\eqref{equ:head-prv-valve-exp} and~\eqref{equ:head-fcv-valve-exp} can be obtained after exponentiating both side of~\eqref{equ:head-prv-valve} and~\eqref{equ:head-fcv-valve}. 
	\begin{linenomath*}
		\begin{subnumcases}{~\label{equ:head-prv-valve-exp}}
		{\hat{h}_{i}} {\hat{h}_{j}^{-1}} [\hat{c}^{\mathrm{W}}]^{-1} [\hat{q}_{ij}^{\mathrm{W}}]^{-1} = 1,\;\mathrm{OPEN} ~\label{equ:head-prv-valve-expa}\\
		\hat{h}_{j}^{-1}  \hat{h}^{\mathrm{W}_\mathrm{set}}= 1,\;\text{ACTIVE} ~\label{equ:head-prv-valve-expb}
		\end{subnumcases}
	\end{linenomath*}
	\begin{linenomath*}
		\begin{subnumcases}{~\label{equ:head-fcv-valve-exp}}
		{\hat{h}_{i}} {\hat{h}_{j}^{-1}} [\hat{c}^{\mathrm{W}}]^{-1} [\hat{q}_{ij}^{\mathrm{W}}]^{-1} = 1,\;\mathrm{OPEN}\label{equ:head-fcv-valve-expa} \\
		{[\hat{q}_{ij}^{\mathrm{W}}]}^{-1} \hat{q}^{\mathrm{W}_\mathrm{set}} = 1,\;\text{ACTIVE} \label{equ:head-fcv-valve-expv}
		\end{subnumcases}
	\end{linenomath*}
	{where $\hat{c}^{\mathrm{W}}(q_{ij}^\mathrm{W})=b^{q_{ij}^\mathrm{W} \left(l_{ij} {|q_{ij}^\mathrm{W}|} - 1\right)}$ in~\eqref{equ:head-prv-valve-exp} and~\eqref{equ:head-fcv-valve-exp}.}
	
	\color{black}
	\subsubsection{Physical constraints}
	For the physical constraints~\eqref{equ:constraints}, the conversion process is similar to Section~\ref{sec:massbalance} since both are linear constraints. After exponentiating~\eqref{equ:constraints},  the GP form becomes
	\begin{linenomath*}
		\begin{subequations}~\label{equ:phsicalConst-exp}
			\begin{align}
			\hat{h}_{i}^{-1}	\hat{h}_{j}^{\mathrm{min}} \leq  1,{\hat{h}_{i}}  \left[ \hat{h}_{j}^{\mathrm{max}}\right]^{-1} &\leq 1 \\
			{\hat{q}_{ij}}^{-1}	{\hat{q}_{ij}}^{\mathrm{min}} \leq  1,{\hat{q}_{ij}}  \left[ {\hat{q}_{ij}}^{\mathrm{max}}\right]^{-1} &\leq 1 \\
			\hat{h}_{ij}^{\mathrm{M}} &\leq 1.~\label{equ:head-gain-exp}
			\end{align}
		\end{subequations} 
	\end{linenomath*}
	\subsubsection{GP modeling of WFP } ~\label{sec:Algorithm}
	After the conversion of all variables and constraints, we can express the converted problem as
	\begin{linenomath*}
		\begin{align}
		\textsc{\textbf{WFP-GP:}}\;   \mathrm{find} \; & \hspace{2em} \hat{\m \xi} \label{equ:WFP-GP} \\
		\mathrm{s.t.}\;&~\eqref{equ:nodes-exp}-\eqref{equ:phsicalConst-exp}. \notag
		\end{align} 
	\end{linenomath*}
	
	Problem~\eqref{equ:WFP-GP} is in standard GP form and can be solved directly by modern GP solvers and even though the WFP-GP~\eqref{equ:WFP-GP} is not convex, as we mentioned, the $\log$ form of this problem is convex{~\cite[Section 2.5]{boyd2007tutorial}}. Starting with an initial guess for the flow rates and relative speeds, the constraints \eqref{equ:nodes-exp}-\eqref{equ:phsicalConst-exp} are approximated at every iteration based on the previous iterations. This process continues until a termination criterion is met. The details are further discussed in Algorithm~\ref{alg:gp-alg}. In the next section, we show how the  WFP-GP~\eqref{equ:WFP-GP}  problem can be formulated using a tractable linear approximation.

	\color{black}
	\subsection{LP modeling derived from GP modeling} \label{sec:LP-GP}
	As we illustrated in the end of Section~\ref{sec:gp-based}, an LP and a GP can be converted to each other, e.g, Eq.~\eqref{equ:GP-linear} can be obtained from~\eqref{equ:LP} via exponent technique, and the inverse operation also holds true meaning that Eq.~\eqref{equ:LP} can be converted back from~\eqref{equ:GP-linear} via the $\log$ operator. Inspired by this idea, we can apply the inverse operation to WFP-GP~\eqref{equ:WFP-GP}, and generate an LP for the WFP.
	After applying the $\log$ function on both sides of GP form of the mass balance equation \eqref{equ:nodes-exp}, the result would be the original linear form of the mass balance equation~\eqref{equ:nodes}. It means the linear modeling is converted back from GP modeling. Similarly, for equations \eqref{equ:head-prv-valve-exp},~\eqref{equ:head-fcv-valve-exp}, and~\eqref{equ:phsicalConst-exp}, the corresponding results are~\eqref{equ:head-prv-valve},~\eqref{equ:head-fcv-valve}, and \eqref{equ:constraints}.

	As for the GP form of energy balance equations for pipes, after applying the $\log$ function with base $b$ on both sides of Eq.~\eqref{equ:head-loss-pipe-exp}, a linear equation can be obtained as
	\begin{linenomath*}
		\begin{equation} ~\label{equ:head-loss-pipe-lp-exp}
		h_{i}-h_{j}-q_{i j}^{\mathrm{P}}=c_{i j}^{\mathrm{P}},
		\end{equation} 
	\end{linenomath*}
	where $c_{i j}^{\mathrm{P}} = \log_b(\hat{c}_{i j}^{\mathrm{P}})$.
	
	Similarly, the linear form of energy equation balance for pumps and valves (GPVs, PRVs, and FCVs) can be expressed as~\eqref{equ:head-flow-pump-lp-exp} and~\eqref{equ:head-flow-gpv-lp-exp} when applying the $\log$ function on both sides of~\eqref{equ:head-flow-pump-exp},~\eqref{equ:head-loss-valve-gpv-exp},~\eqref{equ:head-prv-valve-expa}, and~\eqref{equ:head-fcv-valve-expa}, as
	\begin{linenomath*}
		\begin{align}  
		h_{i}-h_{j}-c_{2}^{\mathrm{M}} q_{ij}^\mathrm{M}={c}_{1}^{\mathrm{M}} \label{equ:head-flow-pump-lp-exp}\\
		h_{i}-h_{j}-q_{i j}^{\mathrm{W}}=c_{i j}^{\mathrm{W}},\label{equ:head-flow-gpv-lp-exp}
		\end{align} 
	\end{linenomath*}
	where $c_{1}^{\mathrm{M}} = \log_b(\hat{c}_{1}^{\mathrm{M}})$, and $c_{i j}^{\mathrm{W}} = \log_b(\hat{c}_{i j}^{\mathrm{W}})$. {Note that \textit{(i)} $\hat{c}_{i j}^{\mathrm{W}}$ varies according to the types of valves, \textit{(ii)} Eq.~\eqref{equ:head-flow-gpv-lp-exp} applies to PRVs and FCVs in open status, and for active status, Eqs.~\eqref{equ:head-prv-valve-active} or~\eqref{equ:head-fcv-valve-active} are used, which are linear as well.} Thus, the nonlinearities from pipes~\eqref{equ:head-flow-pipe}, pumps~\eqref{equ:head-flow-pump}, and valves (GPVs~\eqref{equ:head-flow-valve}, PRVs~\eqref{equ:head-prv-valve}, and FCVs~\eqref{equ:head-fcv-valve}) are approximated by its linear form~\eqref{equ:head-loss-pipe-lp-exp},~\eqref{equ:head-flow-pump-lp-exp}, and~\eqref{equ:head-flow-gpv-lp-exp}. 
	
	{After updating the model in~\eqref{equ:nonlinearMatrix}, the linear matrix representation of WDNs can be written as}
	\begin{linenomath*}
		\begin{equation}~\label{equ:LinearMatrix}
		\hspace{-1em}\includegraphics[width=\linewidth,valign=c]{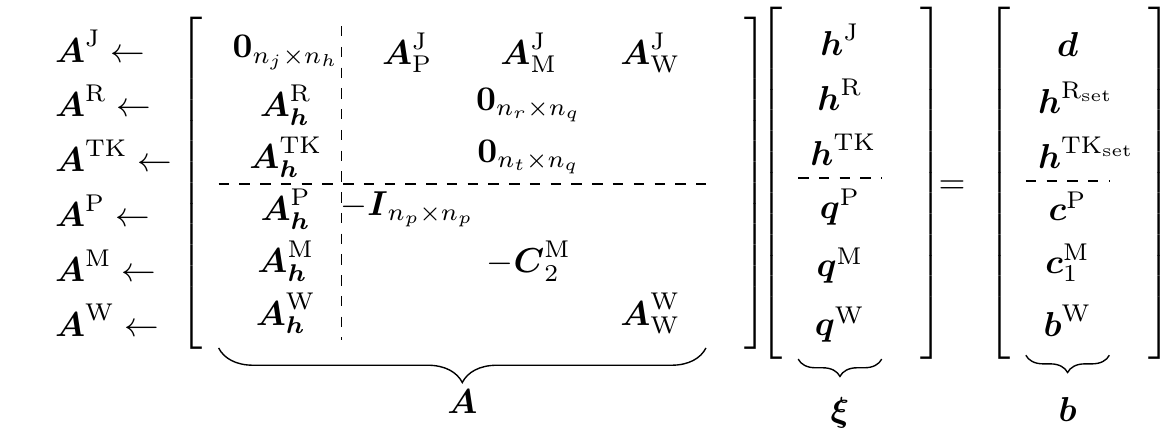},
		\end{equation}
	\end{linenomath*}
	{where ${\m c}^{\mathrm{P}}$, ${\m c}_1^{\mathrm{M}}$,  ${\m C}_2^{\mathrm{M}}$, and ${\m b}^{\mathrm{W}}$ collect the parameters from linear modeling of pipes, pumps, and valves, and note that ${\m c}^{\mathrm{P}} \in \mathbb{R}^{n_{p}}$, ${\m c}_1^{\mathrm{M}}$ and ${\m c}_2^{\mathrm{M}} \in \mathbb{R}^{n_{m}}$, ${\m C}_2^{\mathrm{M}} = \diag({\m c}_2^{\mathrm{M}})  \in \mathbb{R}^{n_m \times n_m}$, and ${\m b}^{\mathrm{W}} \in \mathbb{R}^{n_w}$ including ${\m c}^{\mathrm{W}}$ or valve settings $h^{\mathrm{W}_\mathrm{set}}$, $q^{\mathrm{W}_\mathrm{set}}$. }

	Thus, the LP form of WFP can be expressed as 
	\begin{linenomath*}
		\begin{align}
		\textsc{\textbf{WFP-LP:}}\;  \;\;\;\; \mathrm{find} \;\; \;\;& \hspace{2em}{\m \xi}  \label{equ:WFP-LP} \\
		\mathrm{s.t.}\;\;\;\;& ~\eqref{equ:constraints},~\eqref{equ:LinearMatrix}. \notag
		\end{align}
	\end{linenomath*}
	
	\begin{remark}\label{rm:2}
		$\m A$ in~\eqref{equ:LinearMatrix} is a square matrix, and  $\m A \in \mathbb{R}^{n_{\xi} \times n_{\xi}} $ since  two zero submatrices in $\m A_0$~\eqref{equ:nonlinearMatrix} are replaced with the matrices $-\m I_{n_p \times n_{p}}$ and ${\m C}_2^{\mathrm{M}}$ with the same size.  
	\end{remark}
	\begin{lem}\label{asp:2}
		$\m A$ is invertible. The proof is given in Appendix~\ref{sec:appendixA}.
\end{lem} }

\normalcolor
WFP-LP~\eqref{equ:WFP-LP} derived from GP form can be viewed as a linear approximation of  nonlinear WFP modeling~\eqref{equ:WFP}, and it can be solved with any LP solvers directly. 
Note that all constraints of~\eqref{equ:WFP-LP} are equality constraints except the lower and upper bounds in~\eqref{equ:constraints}, and we can rewrite them in matrix form
\begin{linenomath*}
	\begin{align} ~\label{equ:WFP-LP-matrix}
	&\underbrace{\begin{bmatrix}
		{\m A}^{\mathrm{J}}; & {\m A}^{\mathrm{R}}; &{\m A}^{\mathrm{TK}};  &{\m A}^{\mathrm{P}}; & {\m A}^{\mathrm{M}}; & {\m A}^{\mathrm{W}} 
		\end{bmatrix}}_{\m A} {\m \xi} \\
	&= \underbrace{\begin{bmatrix}
		{\m b}^{\mathrm{J}}; & {\m b}^{\mathrm{R}};  &{\m b}^{\mathrm{TK}};  &{\m b}^{\mathrm{P}}; & {\m b}^{\mathrm{M}}; & {\m b}^{\mathrm{W}}
		\end{bmatrix},}_{\m b} \notag
	\end{align}
\end{linenomath*}
where $\m A$ is a coefficient matrix collecting all submatrices and ${\m A}^{\mathrm{J}}$ is from the mass balance equation~\eqref{equ:nodes}, ${\m A}^{\mathrm{R}} $ and ${\m A}^{\mathrm{TK}}$ collects~\eqref{equ:reservoir} and~\eqref{equ:tank}  for reservoirs and tanks,  ${\m A}^{\mathrm{P}}$ is from the linearized pipe head loss equation~\eqref{equ:head-loss-pipe-lp-exp},  ${\m A}^{\mathrm{M}}$ is linearized equation~\eqref{equ:head-flow-pump-lp-exp} for pump, and ${\m A}^{\mathrm{W}}$ collects all linear equations for valves. The right hand side $\m b$ is a coefficient vector collecting the corresponding coefficients in~\eqref{equ:LinearMatrix}.  An example of $\m A$ and $\m b$ is given in Section~\ref{sec:test}. 

We note that $\m A$ is square, invertible matrix (Lemma~\ref{asp:2}), which implies that an analytical solution can be obtained efficiently large-scale networks using scalable methods for solving linear systems of equations. The bound constraints~\eqref{equ:constraints} are not including in this case, since these constraints are included in design and operation problems to adjust for admissible flows and heads. We will illustrate later in the paper that this approach yields good performance, in comparison with solving a linear program with inequality constraints~\eqref{equ:constraints}.

Next, we provide Algorithm~\ref{alg:gp-alg} for solving the WFP using {WFP-LP}~\eqref{equ:WFP-LP} or its matrix form~\eqref{equ:WFP-LP-matrix}.
Notice that all variables are collected in $\m {{\xi}}$ and the notation $\langle{{\m \xi}}\rangle_{n}$ in Algorithm~\ref{alg:gp-alg} stands for the $n^\mathrm{th}$ iteration value  $\m {{\xi}}$. 
The initial statuses of pumps and valves, head in tanks and reservoirs as well as nodal demands are assumed to be known. For the users familiar with  the EPANET software~\cite{rossman2000epanet}, all the necessary information needed to formulate the WFP-LP~\eqref{equ:WFP-LP} can be seamlessly imported from the ``.inp" source file. 
\begin{algorithm}[t]
	\DontPrintSemicolon
	\KwIn{WDN characteristics/\texttt{.inp} source file, initial guess $\langle{{\m \xi}}\rangle_{0}$, threshold, maxIter, $n_\mathrm{step}$}
	\KwOut{${\m \xi}_{\mathrm{GP-LP}}$}
	Set ${{\m \xi}_{\mathrm{save}}} := \langle{{\m \xi}}\rangle_{0}$, $n=1$, generate $\m A$ and $\m b$ using~\eqref{equ:WFP-LP-matrix}\;
	\While {  $\mathrm{error} \geq \mathrm{threshold}$ \textbf{OR} $n\leq \mathrm{maxIter}$ }{
		Determine the status of each pump and valve and obtain $\langle{{\m c}^{\mathrm{P}}}\rangle_n$, $\langle{\m c_{2}^{\mathrm{M}}}\rangle_n$, and $\langle{{\m c}^{\mathrm{W}}}\rangle_n$ from $\langle{\m {\xi}}\rangle_{n-1}$\;
		Update part of $\m A$, $\m b$, and solve for $\langle{\m {\xi}}\rangle_{n} = {\m A}^{-1} {\m b}$\;
		\If{$\mathrm{mod}(n,{n_\mathrm{step}})=0$}{
			$\m \Delta \m \xi = \langle{\m {\xi}}\rangle_{n}-\langle{\m {\xi}}\rangle_{n-1}$\;
			$\langle{\m {\xi}}\rangle_{n}= \langle{\m {\xi}}\rangle_{n-1} + \m a_n \circ \m \Delta \m \xi $
		}
		Calculate $\mathrm{error} := \mathrm{norm}(\langle{\m {\xi}}\rangle_{n},{\m \xi}_{\mathrm{save}})$\;
		Update ${\m \xi}_{\mathrm{save}}=\langle{\m {\xi}}\rangle_{n}$ and $n=n+1$\;
	}
	Let ${{\m \xi}_{\mathrm{GP-LP}}}=\langle{\m {\xi}}\rangle_{n}$
	\caption{GP-LP method for solving the WFP}
	\label{alg:gp-alg}
\end{algorithm}
{Algorithm~\ref{alg:gp-alg} is initialized with any initial guess $\langle{\m {{\xi}}}\rangle_{0}$ and requires defining the threshold or number iterations for convergence as well as acceleration parameter $\m a_n$. In each iteration, parameters $\langle{{\m c}^{\mathrm{P}}}\rangle_n$ and $\langle{{\m c}^{\mathrm{W}}}\rangle_n$ are updated based on previous iterations, $\langle{\m c_{2}^{\mathrm{M}}}\rangle_n$ is fixed and does not need to be updated, matrices $\m A$ and $\m b$ collecting the above parameters are automatically updated as well}. Notably, for a fixed topology only submatrices ${\m A}^{\mathrm{P}},{\m A}^{\mathrm{M}}$, and parts of ${\m A}^{\mathrm{W}}$ and the corresponding parts in vector $\m b$ require updating, while the rest remain fixed. 
The iteration error is defined as the Euclidean distance between two consecutive iterations. The iterations continue until the error is less than a predefined error threshold ($\mathrm{threshold}$) or the maximum number of iterations ($\mathrm{maxIter}$) is reached, and the final solution is set by ${\m \xi}_{\mathrm{GP-LP}}=\langle{\m {{\xi}}}\rangle_{n}$. 

Steps 5-8 are used to accelerate the convergence of the algorithm. The acceleration parameter $\m a_n$  that can be adjusted dynamically with iterations, e.g. every $n_{\mathrm{step}}$ iterations, and set individually for different elements, where $\m a_n \circ \Delta \m \xi$  represents element-wise product. 
Algorithm~\ref{alg:gp-alg} can also be applied to solve the optimization problem~\eqref{equ:WFP-GP} or~\eqref{equ:WFP-LP} with corresponding GP/LP solver. When dealing with~\eqref{equ:WFP-GP}, the steps remain the same except that the variables are changed into the GP variables $\hat{\m\xi}$. For the same scale problem, the GP solver is usually slower than an LP solver, and the analytical solution is faster than any solver.

\subsection{Convergence of GP-LP iteration}

%

	In this section, we show the convergence of proposed GP method under mild conditions that typically hold in practical WDNs. The  theorem and proof are given first, followed by a discussion on how the initial points and acceleration parameter in Algorithm~\ref{alg:gp-alg} are related with the convergence.
Although Algorithm~\ref{alg:gp-alg} allows checking and updating the status of pumps and valves in each iteration, this section focuses on the typical WFP setting where the statuses are known, and therefore, do not need to be updated in each iteration. In addition, the proof is furnished first for the case where all statuses are open. In this case, GPVs, PRVs, and FCVs have similar modeling of pipes and can be treated directly as pipes with different head loss.  This way, we have a WDN with only pumps and pipes, and valves are included in $n_m$. In what follows, we analyze the iteration $\langle{\m {\xi}}\rangle_{n} = {\m A}^{-1} {\m b}$ of Algorithm~\ref{alg:gp-alg}. 

We set-aside the effect of acceleration parameter at first, that is, we set $\m a_n = \m 0$. Additionally, we combine the equations related with reservoirs and tanks, for example, $\m A^\mathrm{\mathrm{R}}$ and $\m A^\mathrm{\mathrm{TK}}$ are compressed into one row with ${\m h}^{\mathrm{R\_TK}}  = \{	{\m h}^{\mathrm{R}}, {\m h}^{\mathrm{TK}} \}$, $	{\m h}^{\mathrm{set}}  = \{{\m h}^{\mathrm{R}_\mathrm{set}},	{\m h}^{\mathrm{TK}_\mathrm{set}}\}$, and $n_f = n_t + n_r$. Then the matrix in~\eqref{equ:LinearMatrix} takes the form
\begin{figure*}
	\begin{equation}~\label{equ:reogan}
	\hspace{-2.5em}\includegraphics[width=0.88\linewidth,valign=c]{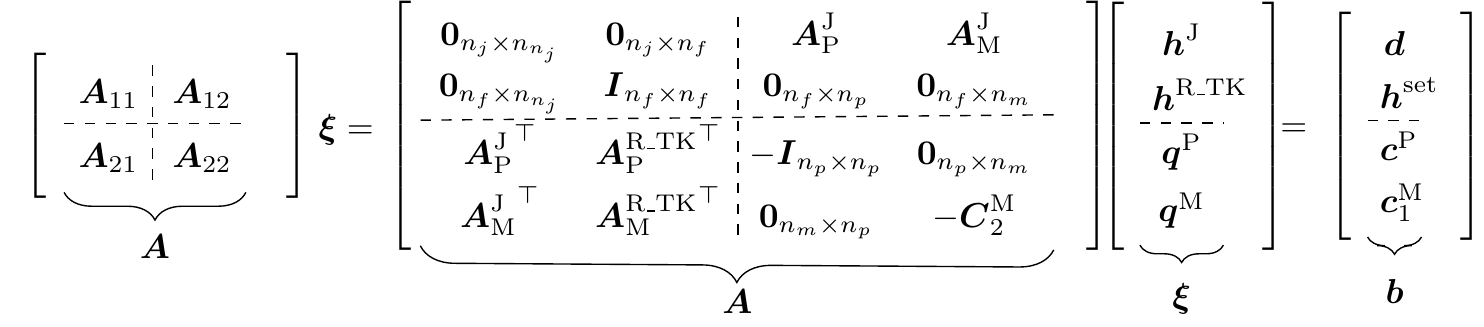}.
	\end{equation}
	\hrulefill
\end{figure*}
Note that $\m A$  is separated into 4 blocks, and $\m A$ changes  if only and if ${\m C}_2^{\mathrm{M}} = \diag({\m c}_2^{\mathrm{M}})$ is updated, where $\m c_2^{\mathrm{M}} = r ({{\m q}^{\mathrm{M}}})^{\nu-1}$ [cf.~\eqref{equ:pumpheadloss}] when speed is fixed as $\m 1$. Similarly, ${\m c}_1^{\mathrm{M}} $ is fixed once the pump curve is set, and $\m b$ changes if only and if $\m c^{\mathrm{P}} = \log_b(\hat{\m c}^{\mathrm{P}}) =\m q^{\mathrm{P}} \circ \left(R {|\m q^{\mathrm{P}}|}^{\mu-1} - \m 1\right)$~\eqref{equ:hat-cp-exp} is updated. In short, $\m A$ only reflects the update of pump flows ${\m q}^{\mathrm{M}} $, and $\m b$ only reflects the updates of pipe flows ${\m q}^{\mathrm{P}}$.

It follows from $\m  \xi = \m A^{-1} \m b$ [cf.~\eqref{equ:LinearMatrix}] that the $n^\mathrm{th}$ iteration is $\langle{\m {\xi}}\rangle_{n} = {\m A}_{n-1}^{-1} {\m b_{n-1}}$, and the $n+1^\mathrm{th}$ iteration  is $\langle{\m {\xi}}\rangle_{n+1} = {\m A}_{n}^{-1} {\m b_{n}}$.	 Now consider two consecutive iterations $n-1$ and $n$, and note that the only changes in $\m A$ and $\m b$ are through the updates $\langle{\m C_2^{\mathrm{M}}}\rangle_{n}  = \langle{\m C_2^{\mathrm{M}}}\rangle_{n-1} + \langle{ \Delta \m C_2^{\mathrm{M}}}\rangle_{n-1}$ and $\langle{\m c^{\mathrm{P}}}\rangle_{n}  = \langle{\m c^{\mathrm{P}}}\rangle_{n-1} +  \langle{\Delta \m c^{\mathrm{P}}}\rangle_{n-1}$. We have $\m A_{n} = \m A_{n-1} +\Delta \m A_{n-1}$ and $\m b_{n} = \m b_{n-1} + \Delta \m b_{n-1}$ where $\Delta \m A_{n-1} = \diag(\m 0, \langle{ \Delta \m C_2^{\mathrm{M}}}\rangle_{n-1} )$ and $\Delta \m b_n = [\m 0 \ \m 0 \  \langle{\Delta \m c^{\mathrm{P}}}\rangle_{n-1} \ \m 0]^\top$. Since we have $n_m$ pumps, we denote the $i^\mathrm{th}$ element in $\Delta \m C_2^{\mathrm{M}}$ (or the parameter for the $i^\mathrm{th}$ pump) as $\Delta c_{2\_i}^{\mathrm{M}}$, and introduce two diagonal matrices $\m U$ and $\m V$ defined as follows:
\begin{equation}
\begin{aligned} 
\m U &= \begin{bmatrix}
\m 0_{(n_h+n_{p})\times (n_h+n_{ p})} &  \m 0\\
\m 0 & \m I_{ n_{m} \times n_{m}}
\end{bmatrix}, \\ 
\m V &= \begin{bmatrix}
\m 0_{(n_h+n_{p})\times (n_h+n_{ p})} &  \m 0\\
\m 0 &\Delta \m C_2^{\mathrm{M}}
\end{bmatrix}.
\end{aligned}
\end{equation}
With the above notation, we have that $\m A_{n} = \m A_{n-1} + \m U \m V$ where $\Delta \m A_{n-1} =  \m U \m V$ is an $n_m$-rank matrix update. Therefore, we have the iterative formula between two consecutive iterations $\langle{\m {\xi}}\rangle_{n+1} = {\m A}_{n}^{-1} (\m b_{n-1} + \Delta \m b_{n-1}) = {\m A}_{n}^{-1} {\m A}_{n-1} \langle{\m {\xi}}\rangle_{n}  + {\m A}_{n}^{-1}  \Delta \m b_{n-1}$, which is written as
\begin{linenomath*}
	\begin{align} \label{equ:iterative}
	\langle{\m {\xi}}\rangle_{n+1} = {\m T}_{n} \langle{\m {\xi}}\rangle_{n}  + { \m e}_{n}  ,
	\end{align}
\end{linenomath*}
where ${\m T}_{n} = {\m A}_{n}^{-1} {\m A}_{n-1}$ and ${\m e}_{n} = {\m A}_{n}^{-1}  \Delta \m b_{n-1}$. 

Introducing an appropriate partition for $\m A_{\mathrm{inv}}=\m A^{-1}$, the vector ${\m e}_{n} = {\m A}_{n}^{-1}  \Delta \m b_{n-1}$ is written as
\begin{linenomath*}
	\begin{equation}~\label{equ:inverse}
	\hspace{-3em}\includegraphics[width=0.92\linewidth,valign=c]{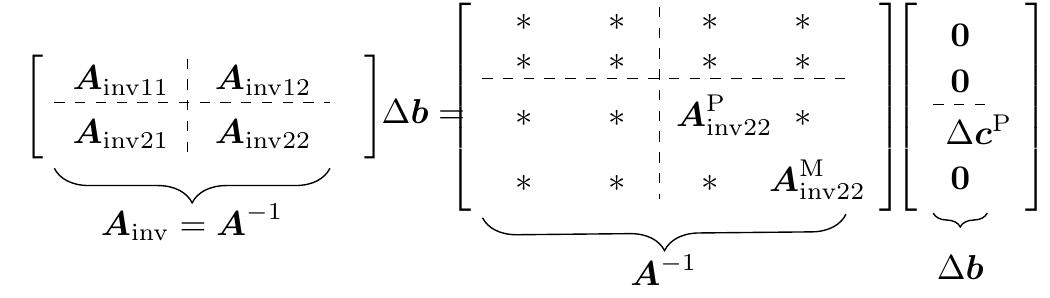},
	\end{equation}
\end{linenomath*}
where $$\m A_{\mathrm{inv}22} = \m A_{22}^{-1} + \m A_{22}^{-1} \m A_{21} [\m A_{11}-\m A_{11}\m A_{22}^{-1}\m A_{21}]^{-1} \m A_{12}\m A_{22}^{-1}$$ according to the block matrix inversion theory~\cite{bernstein2009matrix}. In addition, let $ \m A_{\mathrm{inv}}^{\mathrm{HP,M}}$ denote the entire $(n_h+n_p)\times n_m$ block matrix that sits above  $\m A_{\mathrm{inv}22}^{\mathrm{M}}$. 

Attention is now turned to $\m T_n$. 	According to the Sherman-Morrison-Woodbury formula~\cite{woodbury1950inverting}, we have that
\begin{linenomath*}
	\begin{align*} 
	\m  A_{n}^{-1} &= (\m A_{n-1} +  \m U \m V)^{-1}  \\
	&=\m A_{n-1}^{-1}- {\m A_{n-1}^{-1} \m U { (\m I +\m V \m A_{n-1}^{-1} \m U)^{-1}} \m V} \m A_{n-1}^{-1}.
	\end{align*}
\end{linenomath*}
Introducing the partitions of $\m A_{\mathrm{inv}}$ from~\eqref{equ:inverse},  ${\m T}_{n}$ is written as 
\begin{linenomath*}
	\begin{align} 
	\label{equ:Tn_new}
	{\m T}_{n} \hspace{-2pt} &= \hspace{-2pt} \m I \hspace{-2pt}- \hspace{-2pt} {\m A_{n-1}^{-1} \m U { (\m I \hspace{-2pt}+ \hspace{-2pt}\m V \m A_{n-1}^{-1} \m U)^{-1}} \m V} \\ 
	&= \includegraphics[width=0.7\linewidth,valign=c]{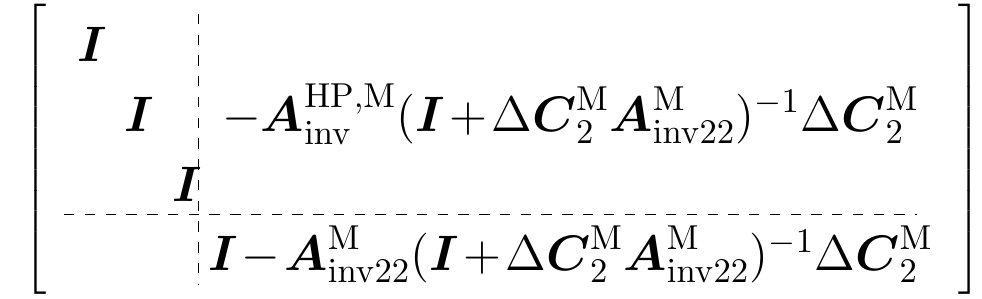} \notag
	\end{align}
\end{linenomath*}
\begin{remark}
	\label{rm:slopes}
	The entries of $\Delta \m C_2^{\mathrm{M}}$ and $\Delta \m c^{\mathrm{P}}$ can be approximated by $\frac{\partial \m C_2^{\mathrm{M}}}{\partial \m q^{\mathrm{M}}} \Delta \m q^{\mathrm{M}}$ and  $\frac{\partial \m c^{\mathrm{P}}}{\partial \m q^{\mathrm{P}}} \Delta \m q^{\mathrm{P}}$.
\end{remark}

{An assumption that facilitates the convergence analysis of~\eqref{equ:iterative} is introduced next. This assumption is expected to be satisfied for typical values of flows in practical WDNs, and was verified numerically for the WDNs of the numerical tests (Section~\ref{sec:test}). 
	
	\begin{assumption}
		\label{as:norm}
		Upon defining the matrix ${\m A_{\mathrm{f}}} \triangleq \diag\left(\mu \m R \circ {|  \m q^{\mathrm{P}} |}^{\mu-1} - \m 1\right)$, it holds that   $\| \m A_{\mathrm{inv}22} {\m A_{\mathrm{f}}} \| <1$, where $\|.\|$ denotes the spectral norm of a matrix.
	\end{assumption}

	The main convergence theorem is stated as follows.

	\begin{thm}
		Under Assumption~\ref{as:norm},  the  GP-LP iteration, or equivalently, the iteration in~\eqref{equ:iterative},  converges. 
	\end{thm}

The proof is given in Appendix~\ref{sec:Convergence}. The condition $\norm{{\m T}^{\mathrm{P}}} = \norm{\langle{\m A_{\mathrm{inv}22}}\rangle_{n} \langle{\m A_{\mathrm{f}}}\rangle_{n}}   < 1$ was valid for all networks we tested and for all $n$. We observe that two factors impact the convergence of GP iteration,  the first one is $\langle{\m A_{\mathrm{inv}22}}\rangle_{n}$ that is mainly decided by the network topology, and the second one is the initialization of the flows that is reflected in matrix $ \langle{\m A_{\mathrm{f}}}\rangle_{n} = \diag\left(\mu \m R {| \langle  \m q^{\mathrm{P}}\rangle_{n-1}|}^{\mu-1} - 1\right)$. When we initialize the flows as zeros, we have that  $\langle{\m A_{\mathrm{f}}}\rangle_{n} = \m I$, which is the worst case. Otherwise, each entry of $\langle{\m A_{\mathrm{f}}}\rangle_{n}$ is in $(-1,0)$, which helps making $\norm{{\m T}^{\mathrm{P}}} < 1$.


{Next, we consider the impact of acceleration parameter $\m a_{n}$ in Algorithm~\ref{alg:gp-alg}. If we increase the previous $ \langle  \Delta \m q^{\mathrm{P}}\rangle_{n-1}$ by $\m a_{n}$, then we have $ \langle{\m A_{\mathrm{f}}}\rangle_{n} = \diag\left(\mu \m R \circ {| \langle  \m q^{\mathrm{P}}\rangle_{n-1} + \m a_n \circ \langle  \Delta \m q^{\mathrm{P}}\rangle_{n-1}|}^{\mu-1} - 1\right)$. In other words,  the $\m a_{n}$ would impact on $\norm{{\m T}^{\mathrm{P}}}$ indirectly via $ \langle{\m A_{\mathrm{f}}}\rangle_{n} $. In order to ensure each entry of $ \langle{\m A_{\mathrm{f}}}\rangle_{n} $ is in $(-1,0)$ and the updated flow is in $[q_{ij}^{\mathrm{min}}, q_{ij}^{\mathrm{max}} ]$, then  corresponding $a_n$ for pipe $ij$ is decided by}
\begin{linenomath*}
	\begin{subequations}~\label{equ:acceleration}
		\begin{align}
		\dfrac{-(\frac{1}{\mu R_{ij}})^{\frac{1}{\mu -1}} - \langle  q_{ij}^{\mathrm{P}}\rangle_{n-1}}{\langle \Delta q_{ij}^{\mathrm{P}}\rangle_{n-1}}<&a_n < \dfrac{(\frac{1}{\mu R_{ij}})^{\frac{1}{\mu -1}} - \langle  q_{ij}^{\mathrm{P}}\rangle_{n-1}}{\langle \Delta q_{ij}^{\mathrm{P}}\rangle_{n-1}}, \label{equ:accelerationA} \\
		\dfrac{ q_{ij}^{\mathrm{min}}- \langle  q_{ij}^{\mathrm{P}}\rangle_{n-1}}{\langle \Delta q_{ij}^{\mathrm{P}}\rangle_{n-1}} <&a_n < \dfrac{ q_{ij}^{\mathrm{max}}- \langle  q_{ij}^{\mathrm{P}}\rangle_{n-1}}{\langle \Delta q_{ij}^{\mathrm{P}}\rangle_{n-1}}.  \label{equ:accelerationB}
		\end{align}
	\end{subequations}
\end{linenomath*}
{From the above $a_n$, we note that \textit{(i)} a pipe with small $R_{ij}$ could be set with large $a_n$. \textit{(ii)} $a_n$ is not only related to the flow rate but also to the change of flow rate in the $(n-1)^\mathrm{th}$ iteration, and when $\langle \Delta q_{ij}^{\mathrm{P}}\rangle_{n-1}$ converges to zero, $a_n$ can theoretically be infinity.  \textit{(iii)} $a_n$  should be limited by~\eqref{equ:accelerationA} and~\eqref{equ:accelerationB} simultaneously. The acceleration parameter $a_n$ that is not limited by the above range can cause the iterations to oscillate or diverge and, hence, a large acceleration parameter can make the iterations less stable.
	
It follows that a proper choice of $\m {a}_n$ needs be chosen for each element and iteration to reach the best performance of convergence. In practice, we find that the iteration starts to oscillate when $\m a_n$ is set to a large value, which makes the overall convergence rate to be slower rather than faster.} We propose adjusting $\m a_n$ in each iteration according to~\eqref{equ:acceleration}, while ensuring that  $\norm{{\m T}^{\mathrm{P}}} < 1$, and reduce $\m a_n$ if the GP iteration start to oscillate.  The above guidelines have been verified by Anytown network in Section~\ref{sec:sensitivity}. Although beyond the scope of this work, the convergence can be optimized by adopting self-adaptation acceleration parameter~\cite{solomon1996accelerating}.

	\section{Case studies}~\label{sec:test}
	Four WDN  examples (3-node, 8-node, Anytown, and C-Town networks) are used to illustrate in detail the applicability of the GP-LP approximation for solving the WFP, and three additional networks are used to test the convergence and simulation times. {The first testcase, that is, the 3-node network, is developed to illustrate the details of GP-LP model.} The second case is a modified 8-node network with a PRV to illustrate that proposed approach is able to handle looped topologies and valves. The C-Town network is used to test the scalability of our approach, and the Anytown network is adopted to discuss the sensitivity analysis. The numerical tests simulated and compared with the help of the EPANET Matlab Toolkit~\cite{eliades2009epanet} on a MacBook Pro with an Intel Core i7 @ 2.2 GHz. No acceleration parameter is used except for the sensitivity analysis of the Anytown network, and different threshold and maximum iterations are set for each network. 
	All the results reported in the next sections are based on solving the LP matrix form of the WFP.  All codes, parameters, tested networks, and results are available on Github~\cite{shenwangGP}.

	\subsection{{Illustrative} 3-node network}
	\begin{figure*}[t]
		\centering
		\includegraphics[width=0.85\linewidth]{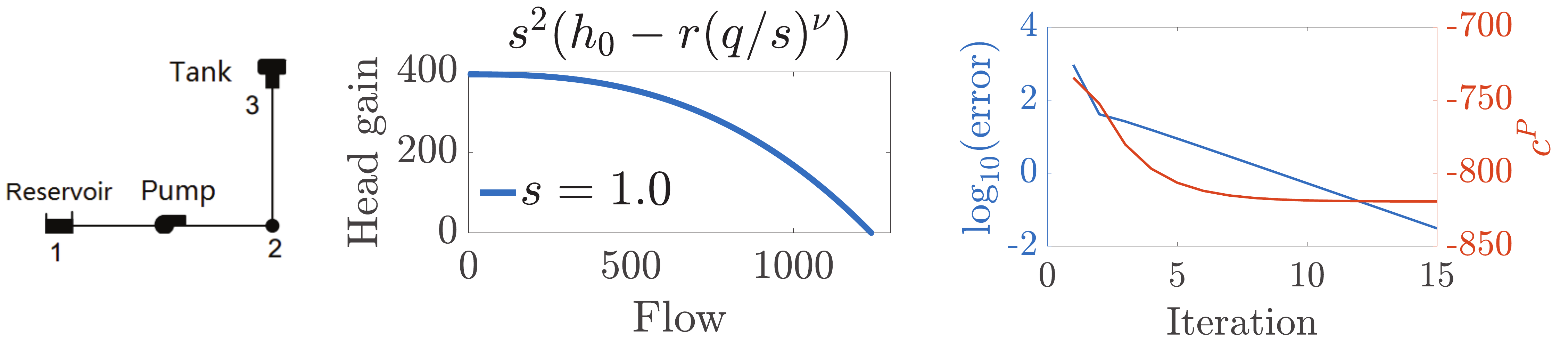}
		\caption{{3-node network (left), variable-speed pump curve (middle), and the value of $\log_{10}(\mathrm{error})$ and ${c}^P$ during iterations (right).}}
		\label{fig:setup}
	\end{figure*}
	In this example, the network is composed of 3 nodes (1 reservoir, 1 tank, and 1 junction with demand) as shown in Fig.~\ref{fig:setup} (left). The corresponding model can be written as~\eqref{eq:3node} in Tab.~\ref{tab:3-node}. {The heads at Reservoir $1$ and Tank $3$ are $h_1^{\mathrm{R}_\mathrm{set}} = 213.4\ \mathrm{m}$, and $ h_3^{\mathrm{TK}_\mathrm{set}} = 276.8\ \mathrm{m}$; the demand at Junction $2$ is $d_2 = 6.3\times 10^{-3}\ \mathrm{m^3/s}$. The curve of variable-speed pump is shown in Fig.~\ref{fig:setup} (middle), the relative speed is known and fixed at $s_{12}= 1$, and the other parameters are $h_0 = 393.7$, $r = 3.8288\times 10^{-6}$, and  $\nu = 2.59$. Note that curve of the pump in Fig.~\ref{fig:setup} (middle) is the negative value of $\Delta h^\mathrm{M}$ defined in~\eqref{equ:head-flow-pump}.} {{{Given the length $L^{\mathrm{P}}=304.8\ \mathrm{m}$, diameter $D^{\mathrm{P}}=0.304\ \mathrm{ft}$, and coefficient $C_{HW} = 100$, the resistance coefficient is $R_{23} = 1.145\times 10^{-5}$. The upper and lower bound constraints are expressed by~\eqref{eq:3node5} to save space.}}} 
	\subsubsection{{Nonlinear modeling of 3-node network}}
	In the problem expressed by Eq.~\eqref{eq:3node}, there are a total of six variables, which can be reduced to three variables because two of them ($h_1$,  and $s_{12}$) are fixed, and $h_3$ represents the water surface elevation because the system is assumed to be operated under steady conditions. The variable $h_2$ can be eliminated by adding~\eqref{eq:3node2} and~\eqref{eq:3node3}. Hence, finding the feasible solution of~\eqref{eq:3node} equals to solving the following nonlinear equations
	\begin{linenomath*}
		\begin{subequations} ~\label{eq:3node-eq}
			\begin{align}
			h_{3} &= h_1^{\mathrm{R}_\mathrm{set}} - R_{23} {q_{23}}\; |q_{23}|^{\mu-1} + (h_0 - r  (q_{12})^\nu ) \label{eq:3node-eq1} \\
			0&= q_{12} - q_{23} - d_2  \label{eq:3node-eq2} \\
			h_{3} &= h_3^{\mathrm{TK}_\mathrm{set}}. \label{eq:3node-eq3}
			\end{align}
		\end{subequations} 
	\end{linenomath*}
	
	The problem represented by~\eqref{eq:3node-eq} which can be visualized in Fig.~\ref{fig:3nodes3D} where each constraint is represented by a corresponding plane. The solution lies in the intersection of surfaces defined by~\eqref{eq:3node-eq1},~\eqref{eq:3node-eq2}, and~\eqref{eq:3node-eq3}. 
	\begin{table*}[t]
		\small
		\caption{3-node network formulation}
		\centering
		\begin{tabular}{cc}
			\hline
			\parbox{6cm}{	
				\vspace{12pt}
				\begingroup
				\setlength\arraycolsep{4pt}
				\begin{subequations} ~\label{eq:3node}
					\begin{align}
					\mathrm{Find} \hspace{33pt}{\boldsymbol \xi}  & \hspace{50pt} \textit{Original form}  \notag  \\
					\mathrm{s.t.} \hspace{5pt} q_{12} - q_{23} &= d_2 \label{eq:3node1} \\
					h_{1} = h_1^{\mathrm{R}_\mathrm{set}},\;h_{3} &= h_3^{\mathrm{TK}_\mathrm{set}}\label{eq:3node4} \\
					h_{2} - h_{3} &= R_{23}\; {q_{23}}\; |q_{23}|^{\mu-1}\label{eq:3node2} \\
					h_{1} - h_{2} &=-(h_0 - r  (q_{12})^\nu ) \label{eq:3node3}\\
					\boldsymbol \xi & \in [  \boldsymbol \xi^{\mathrm{min}},\boldsymbol \xi^{\mathrm{max}} ], \label{eq:3node5}
					\end{align}
				\end{subequations} 
				\endgroup
				\vspace{12pt}
			}& 		\parbox{6cm}{	
				\vspace{-15pt}
				\begingroup
				\setlength\arraycolsep{4pt}
				\begin{subequations}~\label{eq:3node-gp}
					\begin{align}
					\mathrm{Find} \hspace{25pt}  {\hat{\boldsymbol \xi}}  \hspace{70pt} & \hspace{10pt} \textit{GP form}   \notag  \\
					\mathrm{s.t.} \hspace{55pt} \hat{q}_{12}\; \hat{q}_{23}^{-1}\; \hat{d}_2^{-1} &=  1 \label{eq:3node1-gp} \\
					\hat{h}_{1}^{-1} \hat{h}_1^{\mathrm{R}_\mathrm{set}} = 1,\; \hat{h}_{3}^{-1} \hat{h}_3^{\mathrm{TK}_\mathrm{set}} &= 1 \label{eq:3node4-gp} \\
					\hat{h}_{2}\;  \hat{h}_{3}^{-1}\; [\hat{c}^{\mathrm{P}}]^{-1}\; \hat{q}_{23}^{-1}  &= 1      \label{eq:3node2-gp} \\
					\hat{h}_{1}\;  \hat{h}_{2}^{-1}\;  [\hat{c}_1^{\mathrm{M}}]^{-1}\ ({\hat{q}_{12}})^{-c_2^{\mathrm{M}}} &= 1 \label{eq:3node3-gp}\\
					\hat{\boldsymbol \xi} \in [\hat{\boldsymbol \xi}^{\mathrm{min}},\hat{\boldsymbol \xi}^{\mathrm{max}} ] ,& \label{eq:3node5-gp}
					\end{align}
				\end{subequations} 
				\endgroup
				\vspace{-10pt}
			}               \\ \hline
			\parbox{6cm}{	
				\vspace{12pt}
				\begingroup
				\setlength\arraycolsep{4pt}
				\begin{subequations} ~\label{eq:3node-LP}
					\begin{align}
					\mathrm{Find} \hspace{33pt}{\boldsymbol \xi}  & \hspace{80pt} \textit{LP form}  \notag  \\
					\mathrm{s.t.} \hspace{5pt}  q_{12} - q_{23} &= d_2 \label{eq:3nodeLP1} \\
					h_{1} = h_1^{\mathrm{R}_\mathrm{set}},\;h_{3} &= h_3^{\mathrm{TK}_\mathrm{set}}\label{eq:3nodeLP4} \\
					h_{2} - h_{3} &= {c}^{\mathrm{P}} + q_{23}\label{eq:3nodeLP2} \\
					h_{1} - h_{2} &= {c}_1^{\mathrm{M}} +{c_2^{\mathrm{M}}}\ {q}_{12}\label{eq:3nodeLP3}\\
					\boldsymbol \xi & \in [  \boldsymbol \xi^{\mathrm{min}},\boldsymbol \xi^{\mathrm{max}} ], \label{eq:3nodeLP5}
					\end{align}
				\end{subequations} 	
				\endgroup
				\vspace{5pt}
			}
			&  \hspace{-5pt}\parbox{6.2cm}{	
				\vspace{-15pt}
				\begingroup
				\setlength\arraycolsep{1.5pt}
				\begin{align} ~\label{eq:3node-LP_matrix}
				\mathrm{Find} \hspace{33pt}{\boldsymbol \xi}  \hspace{40pt} & \textit{LP matrix form}  \notag  \\
				\underbrace{\begin{bmatrix}
					0 & 0 &0 &1& -\hspace{-1pt}1\\
					0 &1 &0 &0 &0\\
					0 &0& 1& 0& 0 \\
					1 &0& -\hspace{-1pt}1& -\hspace{-1pt}1& 0\\
					-\hspace{-1pt}1& 1& 0& 0& -\hspace{-1pt}c_2^{\mathrm{M}}
					\end{bmatrix}}_{\m A} \hspace{-2pt}\underbrace{\begin{bmatrix}
					h_2\\
					h_1 \\
					h_3\\
					q_{23}\\
					q_{12}
					\end{bmatrix}}_{\m \xi}& \hspace{-2pt}= \hspace{-2pt} \underbrace{\begin{bmatrix}
					d_2\\
					h_1^{\mathrm{R}_\mathrm{set}}\\
					h_3^{\mathrm{TK}_\mathrm{set}}\\
					{c}^{\mathrm{P}}\\
					{c}_1^{\mathrm{M}} 
					\end{bmatrix}}_{\m b}
				\end{align}
				\endgroup
			} \\  \hline
			\label{tab:3-node}
		\end{tabular}
	\end{table*}
	\subsubsection{{GP-LP modeling of 3-node network}}
	The corresponding GP formulation of~\eqref{eq:3node} is listed in~\eqref{eq:3node-gp} in Tab.~\ref{tab:3-node} after applying the technique we introduced in Section~\ref{sec:GPmodeling}.  
	As we mentioned in Section~\ref{sec:LP-GP}, the GP can be transformed to an LP via performing the $\log$ function, as shown in \eqref{eq:3node-LP} in Tab.~\ref{tab:3-node}. Furthermore,~\eqref{eq:3node-LP_matrix} is obtained after rewriting \eqref{eq:3node-LP} in matrix form. We can see that~\eqref{eq:3node-LP} is a linear approximation, but not the same as the first order Taylor approximation. Now this problem can be solved by an LP solver directly. 
	{The parameters we use in Algorithm~\ref{alg:gp-alg} for the 3-node network are selected as: $\mathrm{threshold} = 0.01$ and $\mathrm{maxIter}= 100$.}  Fig.~\ref{fig:setup} (right) shows how the error decreases and ${c}^{\mathrm{P}}$ is updated in each iteration until convergence, which occurs in $n = 15$ iterations. In order to show that our method converges from random initial points, we generate 40 random values for $q_{12}$ and $q_{23}$, and all of them converge to the same final value. {We plot the trajectory of Algorithm~\ref{alg:gp-alg} from random initial guesses to the final solution, represented by the colorful lines in Fig.~\ref{fig:3nodes3D}, where the  marker represents the value of the initial guess, in terms of flow and head, and the color of the line represents the 2-norm distance of current value, where red and blue represent initial guesses farther and closer to the final solution, respectively. We can see that regardless of the the initial value, the solution converges to the final value (blue and small marker). Similar random initializations were performed in the rest of the networks presented in this work demonstrating convergence to the correct solution regardless of the initial guess.} We compare our solution to EPANET simulations and the obtained results are listed in Tab.~\ref{tab:compare}. {The absolute error between $\m \xi_{\mathrm{GP-LP}}$ and $\m \xi_{\mathrm{EPANET}}$ is defined as $\m {\mathrm{AE}} = |\m \xi_{\mathrm{GP-LP}} - \m \xi_{\mathrm{EPANET}}|$,  the corresponding relative error is $\m {\mathrm{RE}}=\frac{\m {\mathrm{AE}}}{| \m \xi_{ \mathrm{EPANET}}|} \times 100\%$, and the Euclidean norm $ {\mathrm{EN}}  = \norm{\m \xi_{\mathrm{GP-LP}} - \m \xi_{\mathrm{EPANET}}}$. The results show that the our approach performs well for this simple tree topology-based network.}

	\begin{table*}
		\small
		\caption{Solution of 3-node network GP-LP versus EPANET.}
		\centering		
		\setlength\tabcolsep{5pt}
		\renewcommand{\arraystretch}{1.5}
		
		\begin{tabular}{cccccc}
			\hline
			\textit{Variables} &$\m \xi_{\mathrm{GP-LP}} $&\hspace{-8pt} $\m \xi_{\mathrm{EPANET}} $ & $\m {\mathrm{RE}}$ & $ {\mathrm{EN}}$
			\\ \hline
			$h_2\ (\mathrm{m})$        & $277.6330$    & $277.6316$            & $0.0005\%$    &     \\
			$q_{12}\ (\mathrm{m^3/s})$    & $5.8186 \times 10^{-2} $  & $5.8186 \times 10^{-2} $          &$0\%$      & $1.4021 \times 10^{-3}$    \\
			$q_{23}\ (\mathrm{m^3/s})$ & $5.1877 \times 10^{-2} $    & $5.1877 \times 10^{-2} $   & $0\%$     &     \\
			\hline
		\end{tabular}
		\label{tab:compare}
	\end{table*}

	\begin{figure}[t]
		\centering
		\includegraphics[width=\linewidth]{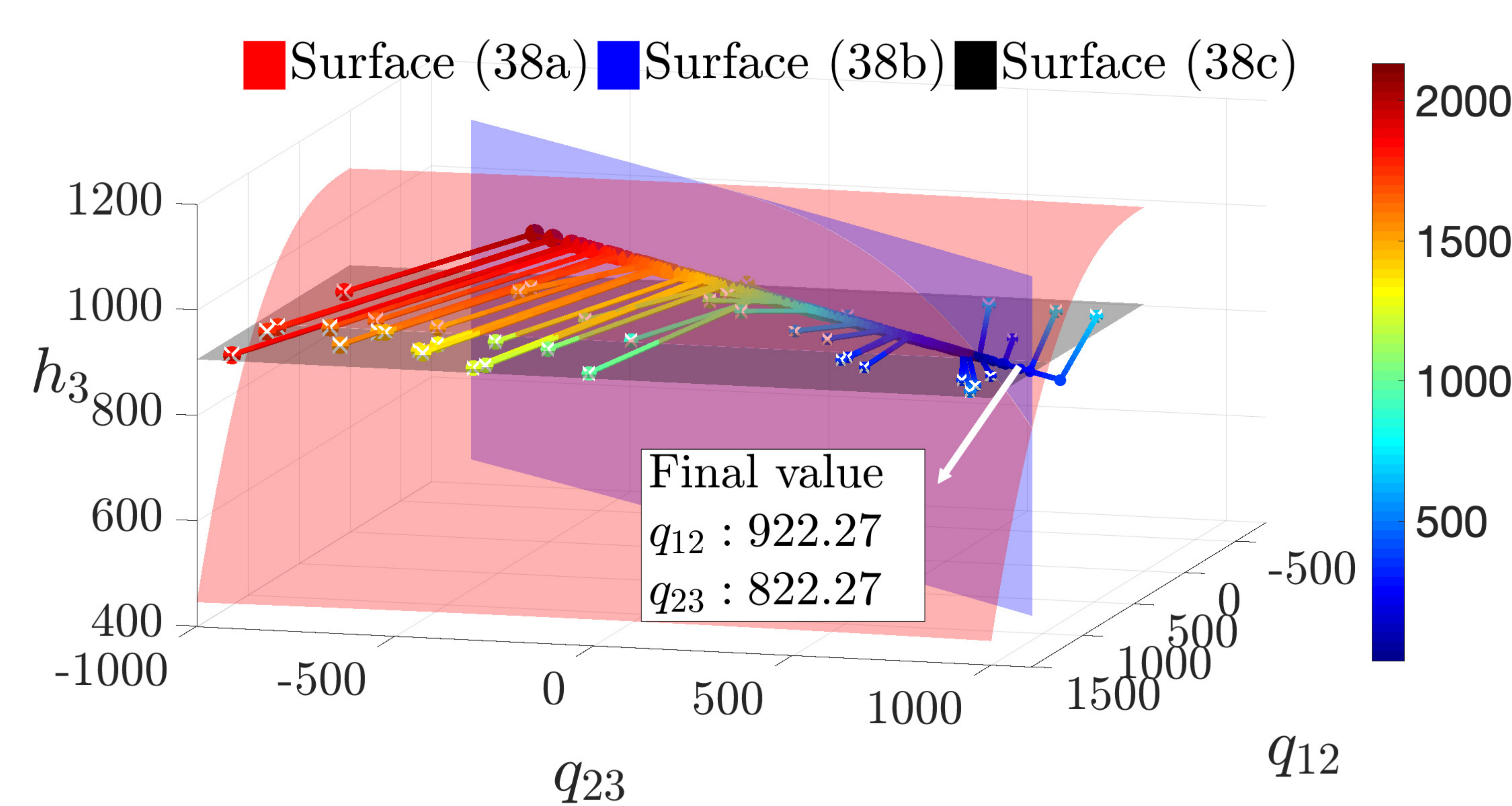} 
		\caption{{Visualization of the iteration process for 40 random initial points}.}
		\label{fig:3nodes3D}
	\end{figure}
	

	
	\subsection{Modified 8-node network}
	The 8-node network is adopted from~\cite[Chapter 2]{rossman2000epanet} and includes a PRV to test our approach with control valves. The modified version includes three more junctions (Nodes 9, 10 and 11) and a PRV between Junctions 3 and 9.  Labels for various components and the topology of modified 8-node network are shown in~Fig.~\ref{fig:complexnetwork}(a). 
	A PRV has two different states corresponding to its working condition. When the PRV is working in ``ACTIVE'' condition, the pressure setting is $P^\mathrm{set} = 45\; \mathrm{m}$. The elevation at downstream side is $E_9 = 190\; \mathrm{m}$. The constraint determined by this PRV from~\eqref{equ:head-prv-valve} is $h_9^\mathrm{W} = 235\; \mathrm{m}$. 
	{The parameters for Algorithm~\ref{alg:gp-alg} are set as: $\mathrm{threshold} = 0.01$ and $\mathrm{maxIter}= 100$.} We test both the ``ACTIVE'' scenario when pressure setting is $45\; \mathrm{m}$ or $100\; \mathrm{m}$ and  the ``OPEN'' scenario.  The error in the three tested scenarios compared to EPANET simulation results is  $\m {\mathrm{EN}} = 0.0067$.   The problem for this test case has  a $23 \times 23$  $\m A$ matrix standing for 23 LP variables and 23 LP constraints when rewritten in the LP matrix form. The difference in number of constraints stems from not considering the upper and lower bound constraints when the network is modeled in the LP matrix form. 
	
	
	\begin{figure}
		\centering
		\includegraphics[width=0.8\linewidth]{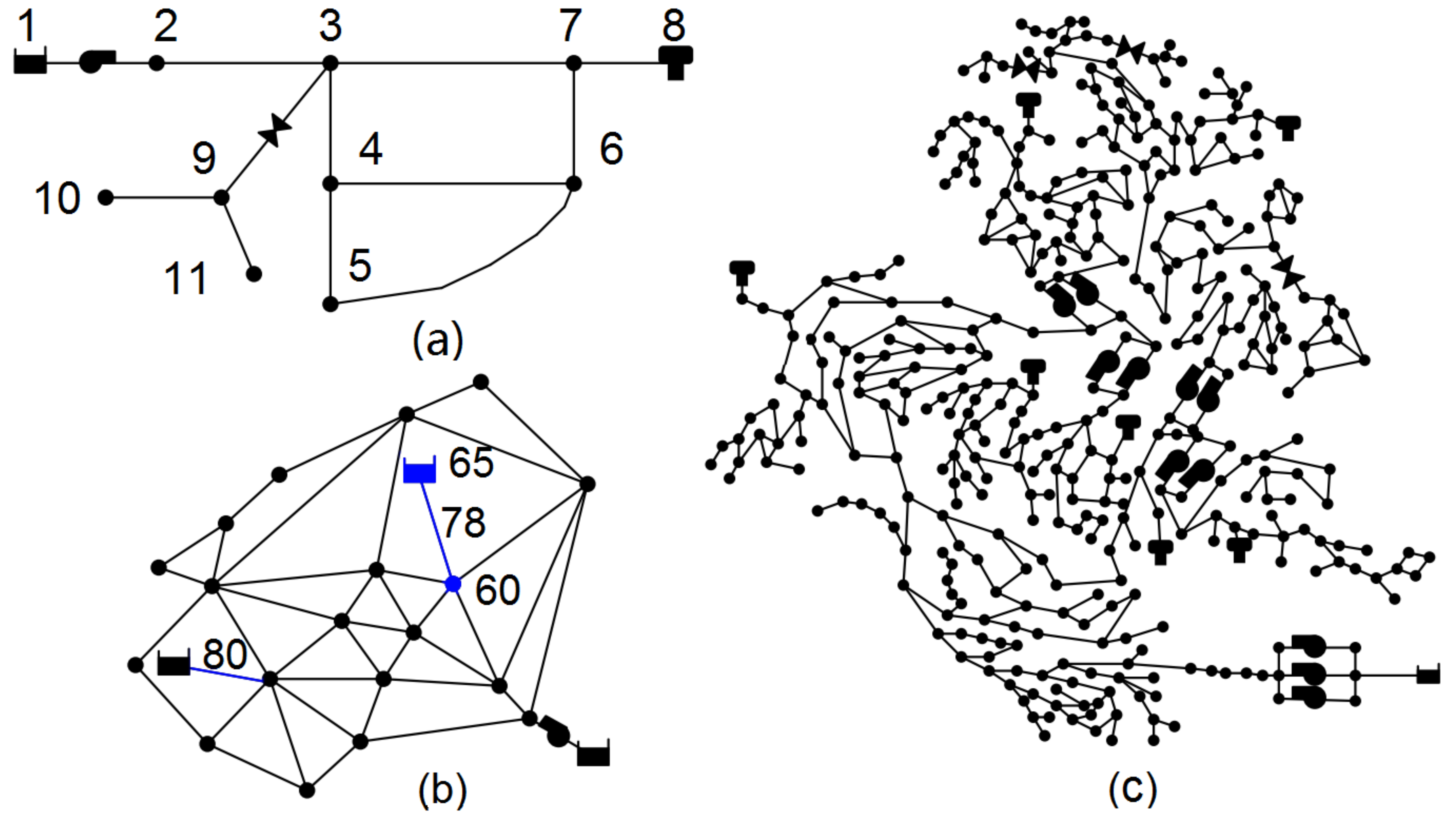}
		\caption{{(a) Modified 8-node network, (b) Anytown, and (c) C-Town network.}}
		\label{fig:complexnetwork}
	\end{figure}
	

	\subsection{C-Town network}

	In order to test the scalability of our proposed GP-LP approach, we test the C-Town network~\cite{ostfeld2012_calibration} \textcolor{black}{that contains 364 junctions, 1 reservoir, 7 tanks, 405 pipes, 11 pumps (three of them are in ``closed" status), and 4 valves  (one of them is closed)} shown in~Fig.~\ref{fig:complexnetwork}\textcolor{black}{(c)}. {The size of $\m A$ is $783\times783$, and the parameters for Algorithm~\ref{alg:gp-alg} are set as: $\mathrm{threshold} = 0.01$ and $\mathrm{maxIter}= 1000$}. 
	
	Fig.~\ref{fig:ctownbar} (left) shows the $\log({\mathrm{EN}})$ with iterations and Fig.~\ref{fig:ctownbar} (right) shows the histogram of absolute errors of individual network components. Note that the convergence criteria is defined as the Euclidean norm between $\m \xi_{\mathrm{GP-LP}} $ and $\m \xi_{\mathrm{EPANET}}$, which summarizes the values of all the components and not individual components. At the final iteration,  $ {\mathrm{EN}} = 1.6969$, which pertains to a $783 \times 1$ vector, thus the error per each variable is small. The histogram in Fig.~\ref{fig:ctownbar} (right) shows that 99\% of absolute errors are within $[0,0.5]$. 
	
	The computational time of solving the LP in matrix form is approximately $5$ $\mathrm{sec}$ with a highly sparse $783\times783$ matrix  $\m A$ in C-Town, where $99.66\%$ of the elements are zeros. The computational time could be further reduced by applying efficient methods for solving $\m A^{-1}$. {Note that this is preliminary work for modeling WDNs, and the majority of computational time involves reading input files, preparing the parameters, and saving temporary results. The code will be optimized in the future work.}
	

	\begin{figure*}
		\centering
		\includegraphics[width=0.8\linewidth]{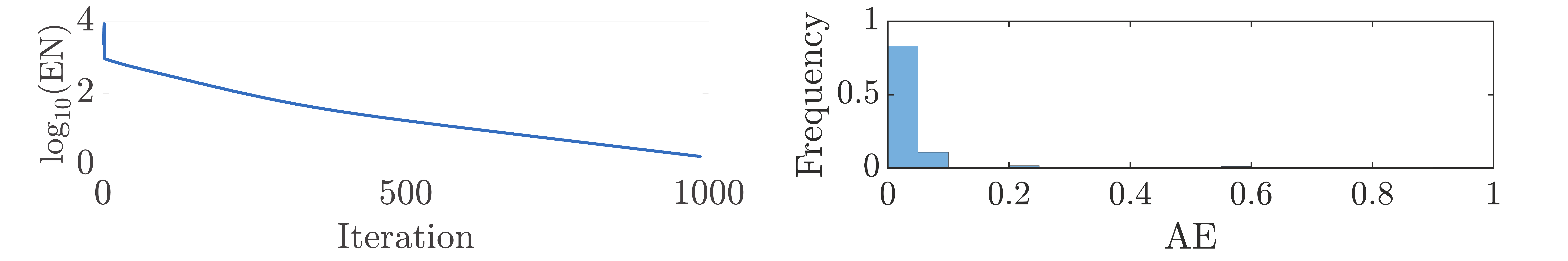}
		\caption{{The $\log_{10}$ of  ${\mathrm{EN}}$ (left) and histogram of  $\m {\mathrm{AE}}$ (right) in C-Town.}}
		\label{fig:ctownbar}
	\end{figure*}

	\subsection{Computational time statistics for tested cases}
	
	We tested the proposed algorithm using additional WDNs that vary in their size and complexity, including PES, NPCL, and OBCL~\cite{eliades2016epanet}. Summary of the main WDNs properties and simulation results are listed in Tab.~\ref{tab:computational-time}. 
	The computational results reported in Tab.~\ref{tab:computational-time} are averaged after testing 3 times for each network with random initial guesses. 
	\begin{table*}[t]
		\small
		\centering
		\caption{Tested networks and the corresponding computational time using Algorithm~\ref{alg:gp-alg}.}
				\setlength\tabcolsep{5pt}
		\renewcommand{\arraystretch}{1.5}
		\label{tab:computational-time}
		\begin{tabular}{c|cccccc}
			\hline
			\textit{Network}            & \textit{8-node} & \textit{Anytown }& \textit{C-Town} &  \textit{PES} & \textit{NPCL} & \textit{OBCL} \\ \hline
			\textit{\makecell{\# of\\components$^*$}}          &   \makecell{\{9,1,1,\\10,1,1\}}   & \makecell{\{19,3,0,\\40,1,0\}}        &  \makecell{\{364,1,7,\\405,11,4\}}       & \makecell{\{68,3,0,\\99,0,0\}}    &  \makecell{\{337,0,2,\\399,0,0\}}    &   \makecell{\{262,1,0,\\288,1,0\}}        \\ \hline
			\textit{\# of variables}          &  23      &   63      &   783     &   170  &  738   &   552         \\
			\textit{\# iterations} &  57     &   {256}     &   729            &   55     &  595    &   87\\
			\textit{Time \textcolor{black}{(sec)}}&  0.0050    &   {0.0553}   &   5.6281         &   0.0346  &  4.3761     &   0.6730\\
			\hline
			\multicolumn{7}{l}{\footnotesize{ \makecell{$^*$\# of components:\{\# Junctions, \# Reservoirs, \# Tanks, \# Pipes, \# Pumps, \# Valves\} }
			}}
		\end{tabular}%
	\end{table*}
	
	\section{Sensitivity analysis}\label{sec:sensitivity}
	To study the sensitivity of our approach and demonstrate how accelerated convergence of Algorithm~\ref{alg:gp-alg} can be achieved, we utilize the Anytown network~\cite{anytown_waslki}. 
	The Anytown network contains 19 junctions, 3 reservoirs, and 40 pipes, as shown in Fig.~\ref{fig:complexnetwork}\textcolor{black}{(b)}. The corresponding LP modeling has has  a   $63 \times 63$ matrix $\m A$. Algorithm~\ref{alg:gp-alg} parameters for Anytown network are set as: $\mathrm{threshold} = 0.01$ and $\mathrm{maxIter}= 5000$. Compared with similar scale networks, the convergence time for Anytown network is relatively slow, requiring 4995 iterations to reach the $\mathrm{threshold}$.  Fig.~\ref{fig:anytownacceleration} shows the change in the error $\log_{10}(\mathrm{EN}) $ with the number of iterations (blue). After analyzing the source data of the Anytown network, we notice that the main errors are caused by the flows through Pipes 78 and 80, connecting the two reservoirs (labeled as blue line segment in~Fig.~\ref{fig:complexnetwork}{(b)}). Here, we only show the analysis for Pipe 78 since both pipes have the same parameters. 
	
	{The first three columns in Tab.~\ref{tab:error_analysis}} show the results after 200 iterations, including the heads at  Reservoir 65 and Junction 60, i.e. $h_{65}$ and $h_{60}$, respectively, the head loss, $\Delta h_{78}^{\mathrm{P}}$, and flow, $q_{78}$, in Pipe 78 connecting Junction 60 with Reservoir 65. Notice that relative error of $h_{60}$ between our final GP-LP solution and EPANET solution is only $0.0092\%$, however, the relative error of $q_{78}$ up to $55.3611\%$. Intuitively, the reason for the error stems from the pipe resistance coefficient, $R$, which is significantly low compared to the rest of the network ($R_{78}=8.1712\times 10^{-7}$). The resistance values of all 40 pipes are plotted in Fig.~\ref{fig:resist}. 
	From the bar plot, we can see that the resistance coefficients for the other pipes in the network are $100$ to $1000$ times greater than for these two pipes. Next, recall the head loss equation~\eqref{equ:head-flow-pipe}, rearranging in terms of the flow we get
	\begin{linenomath*}
		\begin{equation*}
		q_{78}^{\mathrm{P}} = \left( \frac{\Delta h_{78}^{\mathrm{P}}}{R_{78}} \right )^{0.54},
		\end{equation*} 
	\end{linenomath*} 
	where resistance coefficient $R_{78}$ is defined by Hazen-Williams.  Hence, small difference in $\Delta h_{78}^{\mathrm{P}}$ would result in large difference in $q_{78}$ due to the small $R$. For this reason, flows $q_{78}$ and $ q_{80}$ deviate significantly from the EPANET solution. According to the {conservation} of mass~\eqref{equ:nodes}, any variables related to $q_{78}$ and $ q_{80}$ will be affected. 
	\begin{figure}[t]
		\centering
		\includegraphics[width=1.1\linewidth]{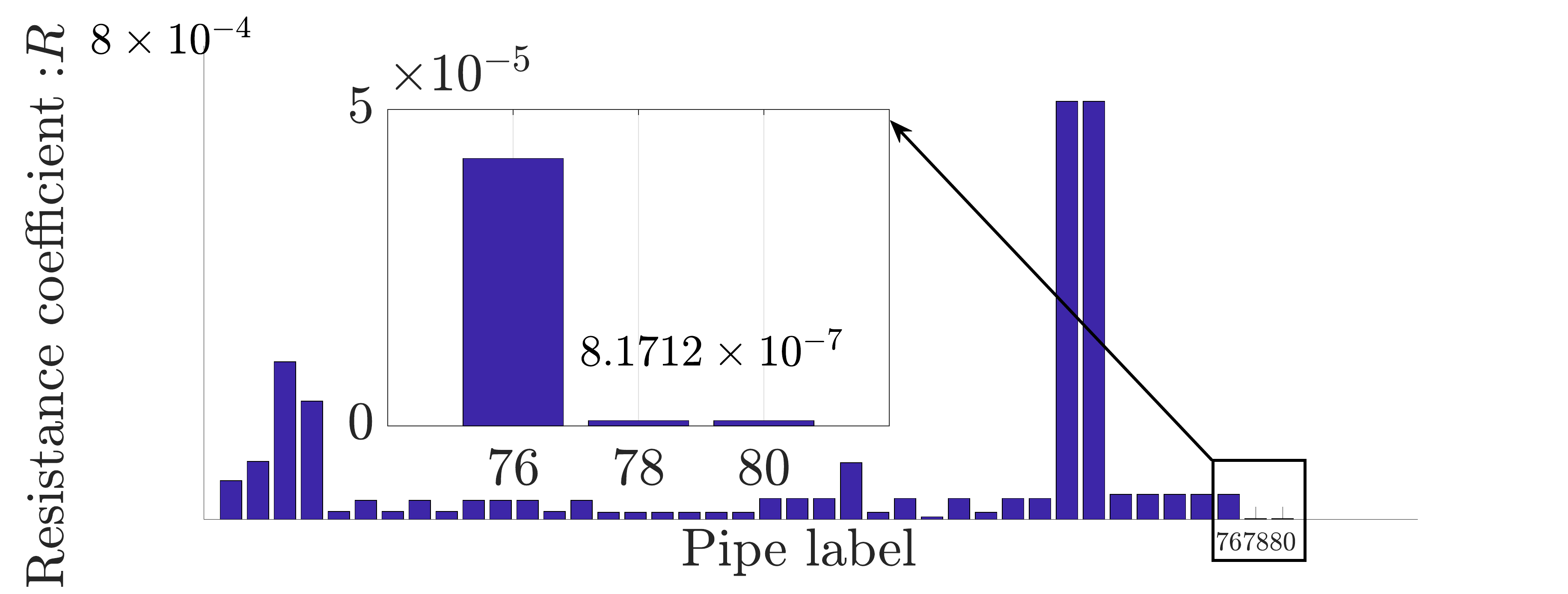}
		\caption{The resistance coefficient $R$ for all 40 pipes in Anytown network (Pipes 76, 78, and 80 are magnified for comparison purposes).}
		\label{fig:resist}
	\end{figure}

	\begin{table*}[t]
		\small
		\centering
		\caption{Sensitivity analysis for the Anytown network.}
		\setlength\tabcolsep{5pt}
		\renewcommand{\arraystretch}{1.5}
		\begin{tabular}{c|ccc|ccc|ccc}
			\hline
			& \multicolumn{3}{c|}{\textit{Base case without $\m a_n$}} & \multicolumn{3}{c|}{\textit{Sensitive analysis}} & \multicolumn{3}{c}{\textit{Base case with $\m a_n$}} \\ \cline{2-10}
			& $\m \xi_{ \mathrm{GP-LP}}$ & $\m \xi_{ \mathrm{EPANET}}$ & $\m {\mathrm{RE}}$ & $\m \xi_{ \mathrm{GP-LP}}$ & $\m \xi_{ \mathrm{EPANET}}$ & $\m {\mathrm{RE}}$ & $\m \xi_{ \mathrm{GP-LP}}$ & $\m \xi_{ \mathrm{EPANET}}$ & $\m {\mathrm{RE}}$ \\ \hline
			$h_{65}$ & 215 & 215 & 0\% & 215 & 215 & 0\%& 215 & 215 & 0\% \\ \hline
			$h_{60}$ &  215.1435& 215.0342 & \textbf{0.0092\%} & 214.4865  & 214.4824  & \textbf{0.0019\%} & 215.0362 & 215.0342 & \textbf{0.0019\%} \\ \hline
			$\Delta h_{78}^{\mathrm{P}}$ &  0.1435 & 0.0342 & ---& 0.5135 &0.5176   &  --- & 0.0362 & 0.0342 & --- \\ \hline
			$q_{78}$ & 139.8081   & 313.1986 & \textbf{55.3611\%}  &-102.8188   & -103.5727  & \textbf{0.7278\%}  & 307.1522 & 313.1986  & \textbf{1.93\%}\\ \hline
		\end{tabular}
		\label{tab:error_analysis}
	\end{table*}

	\begin{figure}[t]
		\centering
		\includegraphics[width=\linewidth]{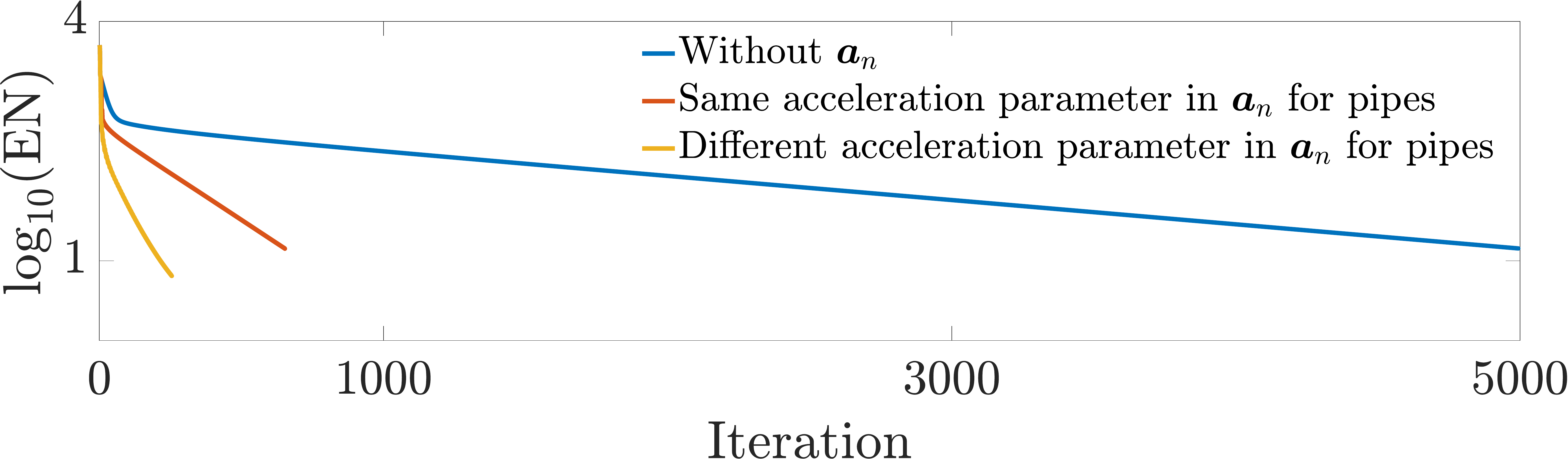}
		\caption{{$\log_{10}(\mathrm{EN}) $ with the number of iterations (4995, 654, and {256}) by adjusting acceleration parameter $\m a_n$.}}
		\label{fig:anytownacceleration}
	\end{figure}

	In order to verify our conclusion that the GP-LP method is sensitive in the resistance coefficients $R$, we increase $R_{78}=9.5873\times 10^{-5}$  by increasing the length, reducing the diameter and the $C_{\mathrm{HW}}$ coefficient of the pipe, and compare  with updated EPANET simulation results.  
	The corresponding results are shown {in the middle three columns of  Tab.~\ref{tab:error_analysis}}, and the relative error of $q_{78}$ after 200 iterations is only $0.7278\%$, where updated parameters are referred as case of sensitivity analysis while the previous parameters are marked as base case. 
	
	To address the sensitivity problem causing the slow convergence, we utilize the acceleration parameter $\m a_n$ in Algorithm~\ref{alg:gp-alg}. We test two scenarios: \textit{(i)} all acceleration parameters are set to 20 and \textit{(ii)} adjusting the acceleration parameters in each iteration according to~\eqref{equ:acceleration}. Fig.~\ref{fig:anytownacceleration} shows the $\log_{10}(\mathrm{EN}) $ with the number of iterations without an acceleration parameter (blue), with the same acceleration parameter for all pipes (red), and adjusting the acceleration parameter for Pipes 78 and 80 (yellow).  Comparing with the base case without $\m a_n$ which takes 4995 iterations to converge, using same acceleration parameters converges withing 654 iterations, and with different acceleration parameters converges within {256} iterations. The latter solution is presented in the last three columns in Tab.~\ref{tab:error_analysis} and the maximum relative error caused by $q_{78}$ is only $1.93\%$, which  improved significantly compared with the base case.
	
	
	\section{Using GP-NET and  future extensions}~\label{sec:future}
	\vspace{-0.75cm}
	\subsection{Using the proposed WFP solver: GP-NET}
	All the codes, tested networks, and results are available on~\cite{shenwangGP}. From our formulations, the code of solving WFP can be divided into three categories. The first category is using WFP-GP~\eqref{equ:WFP-GP} formulation and can be solved by general GP solvers \cite{mosek, cvx}. The second category is formulated based on WFP-LP~\eqref{equ:WFP-LP}  and can be solved is using LP solvers. The third is the matrix form of LP~\eqref{equ:WFP-LP-matrix}, which needs no solvers, and the analytical solution can be obtained for each iteration. Technically, WFP-GP and WFP-LP are more general, and could be adapted to any other problems such as control and state estimation problems~\cite{wang2019geometric} in WDNs by simply modifying the objective function or adding constraints, while the analytical solution provided by matrix form of LP is faster and  can only be applied to solving the WFP. We note that not all nonlinear functions can be converted into GP form so far, and both WFP-GP and WFP-LP are meaningful since some objective functions are easier to formulate as GP or vice versa. Since all formulations are based on or derived from GP,  we name the solver GP-NET. In order to use GP-NET properly, the reader is referred to the Readme.md file on Github~\cite{shenwangGP}.
	
	\subsection{Pressure driven modeling}~\label{sec:PDM}
	Pressure driven demand can be seamlessly integrated into the GP modeling approach.  The pressure driven demand function can be formulated as~\cite{Giustolisi2008,Giustolisi2011}
	\begin{linenomath*}
		\begin{subnumcases}	{\hspace{-2em} d_i^{\mathrm{act}}=~\label{equ:pdd1}}
		d_i^{\mathrm{dsgn}} & $h_i \geq h_i^{\mathrm{ser}}$ \label{equ:norm}
		\\
		d_i^{\mathrm{dsgn}} \left (\frac{h_i - h_i^{\mathrm{min}}}{h_i - h_i^{\mathrm{min}}}\right )^{\gamma} & $h_i^{\mathrm{min}} < h_i < h_i^{\mathrm{ser}}$ \label{equ:reduced}
		\\
		0 & $h_i \leq h_i^{\mathrm{min}}$ \label{equ:fail},
		\end{subnumcases}
	\end{linenomath*}
	where $d_i^{\mathrm{act}}$ is the actual supplied demand, $d_i^{\mathrm{dsgn}}$ is the desired demand,  $h_i^\mathrm{ser}$ and $h_i^{\mathrm{min}}$ are the service and the minimum heads, respectively, and $\gamma$ is typically equal to 0.5~\cite{wagner1988water}. 
	
	
	Given the pressure driven demand model, the mass balance equation~\eqref{equ:nodes} can be rewritten as
	\begin{linenomath*}
		\begin{equation}~\label{equ:pdd}
		\sum_{j \in \mathcal{N}_{i}^{\mathrm{in}}} q_{j i}-\sum_{j \in \mathcal{N}_{i}^{\mathrm{out}}} q_{i j}= d_i^{\mathrm{act}}.
		\end{equation}
	\end{linenomath*}
	The pressure driven demand is a function of $h_i$ in three different regions. The shape of the demand function for $h_i \geq h_i^{\mathrm{ser}}$ and $h_i \leq h_i^{\mathrm{min}}$ is similar to the original demand driven formulation~\eqref{equ:nodes}. For $h_i^{\mathrm{min}} < h_i < h_i^{\mathrm{ser}}$,  the demand in~\eqref{equ:reduced} has the similar form of the head loss model $h_{i}-h_{j}=R_{i j} \left|q_{i j}^\mathrm{P}\right|^{\mu}$ in~\eqref{equ:head-flow-pipe}. In fact, it is easier than the head loss model~\eqref{equ:head-flow-pipe} because~\eqref{equ:reduced} does not have the absolute sign, hence, the same trick we introduced to deal with the head loss model is also applicable. The GP form of~\eqref{equ:pdd} can be expressed as
	\begin{linenomath*}
		\begin{equation*}
		\prod_{j \in \mathcal{N}_{i}^{\mathrm{in}}} \hat{q}_{ji}^{-1} \prod_{j \in \mathcal{N}_{i}^{\mathrm{out}}} \hat{q}_{i j}^{-1} \hat{h}_{i}^{-1}\left[\hat{c}^{\mathrm{J}}\right]^{-1}=1,
		\end{equation*}
	\end{linenomath*}
	where $\hat{c}^{\mathrm{J}} = b^{h_i-d_i^{\mathrm{act}}}$ is a parameter similar as $\hat{c}^{\mathrm{P}}$ in the GP form of head loss. The corresponding LP form is written as
	\begin{linenomath*}
		\begin{equation*}
		\sum_{j \in \mathcal{N}_{i}^{\mathrm{in}}} q_{j i}-\sum_{j \in \mathcal{N}_{i}^{\mathrm{out}}} q_{i j} - h_i= {c}^{\mathrm{J}}.
		\end{equation*}
	\end{linenomath*} 
	Note that ${c}^{\mathrm{J}}$ is updated similarly as the updating process in~\eqref{equ:head-prv-valve-exp} and~\eqref{equ:head-fcv-valve-exp}, by checking the value of $h_i$ in the previous iteration and selecting the appropriate function  $d_i^{\mathrm{act}}$ to update the new value of ${c}^{\mathrm{J}}$. Nodal leakage can be modeled similarly to~\eqref{equ:pdd1}, having two cases for negative and positive nodal pressure head, and by superimposing leakage from contributing pipes proportionally to the pressure at the incident nodes~\cite{Giustolisi2008}.
	
	\subsection{Optimal control} ~\label{sec:OPC}
	The WFP can be integrated in optimization problems for different applications. 
	For example, consider optimal tank and pump control in WDNs, in which the objective of the operator is to minimize the deviation of the water levels in tanks $\m h^\mathrm{TK}(k)$ from a target value $\m h^{\mathrm{TK}_{\mathrm{set}}}$ or enforce smooth operation by minimizing the variability of pump operations $\Delta \m q^\mathrm{M}(k)$. These can be added as objective functions to form a WFP-constrained optimization problem formulated in WFP-GP~\eqref{equ:WFP-GP} or in WFP-LP~\eqref{equ:WFP-LP}. Specifically,
	\begin{linenomath*}
		\begin{subequations}
			\begin{align}
			\hspace{-5pt}&\Gamma_1(k) =  \left(\m h^{\mathrm{TK}}(k)- \m h^{\mathrm{\mathrm{TK}}_{\mathrm{set}}}\right)^{\top}\hspace{-3pt}\left(\m h^{\mathrm{TK}}(k)  - \m h^{\mathrm{TK}_{\mathrm{set}}} \right)~\label{equ:SafetyWater} \\
			\hspace{-5pt}&\Gamma_2(k) ={ \Delta \m q^\mathrm{M}(k)}^{\top} { \Delta \m q^\mathrm{M}(k)}, ~\label{equ:SmoothinControl}
			\end{align}
		\end{subequations}
	\end{linenomath*}
	where $\Gamma_1(\cdot)$ promotes maintaining the targeted water storage set by the operator;  $\m h^{\mathrm{TK}}$ collects the head in tanks, and $\m h^{\mathrm{TK}_{\mathrm{set}}}$ is a vector collecting the target head levels of tanks; $\Gamma_2(\cdot)$ promotes the smoothness of control actions through $\Delta \m q^\mathrm{M}(k) = \m q^\mathrm{M}(k) - \m q^\mathrm{M}(k-1)$ by minimizing the variability in the flow rate changes of controllable components from time $k-1$ to $k$.

	We can convert the above objective functions using proposed GP tricks: \textit{(1) Conversion of $\Gamma_{1}$: }
	The objective	$\Gamma_{1}$ promotes $\m h^{\mathrm{TK}}$ to be close to $\m h^{\mathrm{TK}_{\mathrm{set}}}$. Hence, we introduce a new auxiliary variable $\hat{\m z}(k) \triangleq b^{\m h^{\mathrm{TK}_{\mathrm{set}}} - \m h^{\mathrm{TK}}(k)}$ which will be close to $\m 1$ when tank water levels are close to the target levels. Using the epigraph form, the original objective function $\Gamma_1$  is replaced with $\hat{\Gamma}_1(\hat{\m z}(k)) = \prod_{i = 1}^{n_t}  \hat{z}_i(k)$ and the following constraints are added $\hat{z}_i(k) = \hat{h}_i^{\mathrm{TK}_{\mathrm{set}}} [\hat{ h}_i^{\mathrm{TK}}(k)]^{-1}$ and $\hat{z}_i(k) \geq 1$, where $ \hat{h}_i^{\mathrm{TK}_{\mathrm{set}}}=b^{h_i^{\mathrm{TK}_{\mathrm{set}}}}$ and $ \hat{h}^{\mathrm{TK}}(k)=b^{h^{\mathrm{TK}}(k)}$.  \textit{(2) Conversion of $\Gamma_{2}$: }
	Using the epigraph form, the original objective function $\Gamma_2$ can be expressed as a new objective $\hat{\Gamma}_2(\hat{\m p}(k)) = \prod_{i = 1}^{n_m}  \hat{p}_i(k)^{\Delta q^\mathrm{M}_i(k)}$ with additional constraints given as $\hat{ p}_i(k) = \hat{q}^\mathrm{M}_i(k) [\hat{q}^\mathrm{M}_i(k-1)]^{-1} $ and $\hat{\Gamma}_2(\hat{\m p}(k)) \geq \beta$ where parameter $\beta$ stands for the extent of smoothness of the objective function and to prevent $\hat{\Gamma}_2(\hat{\m p}(k))$ from reducing to 0. For more details, the reader is referred to our recent paper~\cite{tcns} where we thoroughly investigate optimal pump and valve control jointly with the presented GP-based methods in the present paper.
	
	\section{Paper Summary and Future Work}\label{sec:conclusions}
	In this paper, a new derivative-free, linear approximation method is proposed for solving the water flow problem. The proposed approach transforms the variables and constraints in WFP, which, ultimately, reduces to an LP that can be  solved analytically or by linear solvers.   Case studies demonstrate the performance in terms of accuracy and convergence rates. The proposed approach considers looped and branched network topologies, flow directions, and various valve types, it is scalable to large water networks. Under mild conditions, we show that the proposed linear approximation converges and provide guidelines to achieve convergence speedups.  Additionally, we demonstrate future extensions to include pressure driven demand and leak modeling as well as integrating the WFP in network control and state estimation problems.
	
	 The modeling approach proposed in this work can be transferred to other infrastructure systems, such as natural gas infrastructure, in which the governing equations can be modeled similarly to water networks, where gas flow and pressures correspond to water flow and heads, compressor and regulator stations correspond to pumps and control valves that increase and regulate pressures, respectively~\cite{6954481, gas2012, gas2015}. Future work will explore possible extensions to other infrastructure systems as well as further improving the computational performance of the proposed approach.

	
	\normalcolor

	\appendices
	
	\section{Proving that \texorpdfstring{$\m A$}{TEXT}  is invertible} \label{sec:appendixA} 
	\begin{proof}
		According to Assumption~\ref{asp:1}, we know that $\m A_0$ is invertible, i.e. a row of all zeros does not exist after Gaussian elimination. Hence, the row submatrices  in $\m A_0$ are linearly independent with each other, which means $\m A_{\m h}^\mathrm{R}$, $\m A_{\m h}^\mathrm{TK}$, $\m A_{\m h}^\mathrm{P}$, and $\m A_{\m h}^\mathrm{M}$ are linearly independent with each other.
		According to Remark~\ref{rm:2}, row submatrices ${\m A}^{\mathrm{J}}$,
		${\m A}^{\mathrm{R}}$,
		${\m A}^{\mathrm{TK}}$,
		and ${\m A}^{\mathrm{W}}$ in $\m A$ are also linearly independent with each other. Thus, in order to prove $\m A$ is invertible, we only need to prove ${\m A}^{\mathrm{P}}$ and ${\m A}^{\mathrm{M}}$ are linearly independent with the other row submatrices in $\m A$. We will prove that ${\m A}^{\mathrm{P}}$ (corresponding to the pipes) is  linearly independent of the rest of the submatrices.
		
		First, we can see that ${\m A}^{\mathrm{P}} \in \mbb{R}^{n_p \times n_p}$ itself is linearly independent because it contains a $\m I_{n_p \times n_{p}}$. Second, ${\m A}^{\mathrm{P}}$ is linearly independent of ${\m A}^{\mathrm{R}}$, ${\m A}^{\mathrm{TK}}$, ${\m A}^{\mathrm{M}}$, and ${\m A}^{\mathrm{W}}$ because identity matrix $\m I_{n_p \times n_{p}}$ can not be eliminated with zero rows using Gaussian elimination. {Third, it is clear that ${\m A}^{\mathrm{P}}$ is linearly independent with ${\m A}^{\mathrm{J}}$ because each row in ${\m A}^{\mathrm{P}}$ collecting~\eqref{equ:head-loss-pipe-lp-exp} includes linear combination of heads and flows, while each row in ${\m A}^{\mathrm{J}}$ collecting~\eqref{equ:nodes} includes linear combination of flows.}
		Similarly, we can prove that ${\m A}^{\mathrm{M}}$ is  linearly independent of  the rest of submatrices. Hence, $\m A$ is invertible.
	\end{proof}
	
	\section{Convergence proof of the GP-LP iteration}~\label{sec:Convergence}
			\vspace{-0.5cm}
		\begin{proof}
			The convergence of vector $\m q^{\mathrm{P}}$ is shown first. {With Remark~\ref{rm:slopes}, }for the $i^\mathrm{th}$ pump, $\Delta c_{2\_i}^{\mathrm{M}}$ from $\Delta \m C_2^{\mathrm{M}}$ is obtained from the slope of $c_{2\_i}^{\mathrm{M}}$ times the changes in  $q_i^{\mathrm{M}}$, that is, $\Delta c_{2\_i}^{\mathrm{M}} = r (q_i^{\mathrm{M}})^{\nu-2} \Delta q_i^{\mathrm{M}} =  c_{2\_i}^{\mathrm{M}}\frac{\Delta q_i^{\mathrm{M}}}{ q_i^{\mathrm{M}}}$. The typical value of $c_{2\_i}^{\mathrm{M}}$ is small due to the fact that parameter $r$ of the pump curve is very small and also renders  $\Delta c_{2\_i}^{\mathrm{M}}$ even smaller than $c_{2\_i}^{\mathrm{M}}$. The typical value of $\Delta c_{2\_i}^{\mathrm{M}}$ is thus small enough ($10^{-5}\ \mathrm{m}$) in practice, which makes the diagonal elements in $\Delta \m C_2^{\mathrm{M}}$ very small. It follows that the block of $\m T_n$ that depends on $\Delta \m C_2^{\mathrm{M}}$ [cf.~\eqref{equ:Tn_new}] becomes negligible. 
			
			Combining the latter with~\eqref{equ:Tn_new} and~\eqref{equ:inverse} , the iteration for the components of $\m q^{\mathrm{P}}$ in~\eqref{equ:iterative} takes the following form:
			\begin{linenomath*}
				\begin{align}  \label{equ:flowq}
				\langle{\m q^{\mathrm{P}}}\rangle_{n+1} = \langle{\m q^{\mathrm{P}}}\rangle_{n}  + \langle{\m A_{\mathrm{inv}22}}\rangle_{n}  \langle{\Delta \m c^{\mathrm{P}}}\rangle_{n-1}.
				\end{align} 
			\end{linenomath*}
			{With Remark~\ref{rm:slopes}}, we have that $ \langle{\Delta \m c^{\mathrm{P}}}\rangle_{n-1} =  \diag\left(\mu \m R \circ {| \langle \m q^{\mathrm{P}}\rangle_{n-1}|}^{\mu-1} - 1\right)\langle \Delta \m q^{\mathrm{P}}\rangle_{n-1}$. Let $ \langle{\m A_{\mathrm{f}}}\rangle_{n} \triangleq \diag\left(\mu \m R \circ {| \langle  \m q^{\mathrm{P}}\rangle_{n-1}|}^{\mu-1} - 1\right)$, and therefore,~\eqref{equ:flowq} becomes
			\begin{linenomath*}
				\begin{align}  \label{equ:flowq1}
				\langle{\Delta \m q^{\mathrm{P}}}\rangle_{n} =  \langle{\m A_{\mathrm{inv}22}}\rangle_{n} \langle{\m A_{\mathrm{f}}}\rangle_{n} \langle \Delta \m q^{\mathrm{P}}\rangle_{n-1} \triangleq {\m T}_{n}^{\mathrm{P}} \langle \Delta \m q^{\mathrm{P}}\rangle_{n-1},
				\end{align}
			\end{linenomath*}
			where ${\m T}_{n}^{\mathrm{P}} =  \langle{\m A_{\mathrm{inv}22}}\rangle_{n} \langle{\m A_{\mathrm{f}}}\rangle_{n} $. Note that each diagonal entry of $\m A_{\mathrm{f}}$ is in $\left[-1,0\right)$,  indeed, in order for the entries of  $\m A_{\mathrm{f}}$ to be outside of the interval $\left[-1,0\right)$, it would be required that  $|q_{ij}^{\mathrm{P}}| \geq \left(\dfrac{1}{\mu R_{ij}}\right)^{1/(\mu-1)}$. {For the typical values of $\mu = 1.852$ and  $R_{ij}=1\times10^{-5}$ (unitless) using Hazen-Williams model,  the latter condition implies $|q_{ij}^{\mathrm{P}}| \geq 24\ \mathrm{CMS}$,} which clearly cannot hold in practical WDNs.  
			Invoking Assumption~\ref{as:norm}, it follows from~\eqref{equ:flowq1} that $\langle{\Delta \m q^{\mathrm{P}}}\rangle_{n} \rightarrow 0$, which implies that $\langle{ \m q^{\mathrm{P}}}\rangle_{n}$ converges.
			
			Attention is now turned to the flows through pumps.  From the conservation of mass~\eqref{equ:nodes}, we know that  $\m q^{\mathrm{M}}$ can always be expressed as the linear combination of $\m q^{\mathrm{P}}$ and demand $\m d$, which is fixed. Hence, $\m q^{\mathrm{M}}$ converges when $\m q^{\mathrm{P}}$ converges. Next we prove the convergence of the $\m h$ components of $\m \xi$.
			
			From~\eqref{equ:reogan} or original~\eqref{equ:nonlinearMatrix}, we can obtain ${\m h}^{\mathrm{R\_TK}} = {\m h}^{\mathrm{set}}$ and 
			\begin{linenomath*}
				\begin{align*}
				\begin{bmatrix}
				{\m A^\mathrm{J}_\mathrm{P}}^\top\\
				{\m A^\mathrm{J}_\mathrm{M}}^\top
				\end{bmatrix} \m h^\mathrm{J} = 
				\begin{bmatrix}
				\Delta \m h^\mathrm{P}(\m q^\mathrm{P})\\
				\Delta \m h^\mathrm{M}(\m q^\mathrm{M})
				\end{bmatrix} - 			 \begin{bmatrix}
				{\m A^\mathrm{R\_{TK}}_{\mathrm{P}}}^\top\\
				{\m A^\mathrm{R\_{TK}}_{\mathrm{M}}}^\top
				\end{bmatrix} {\m h}^{\mathrm{R\_TK}},
				\end{align*}
			\end{linenomath*}
			
			where $\Delta \m h^\mathrm{P}(\m q^\mathrm{P})$ and $\Delta \m h^\mathrm{P}(\m q^\mathrm{M})$ converge when $\m q^{\mathrm{P}}$ and $\m q^{\mathrm{M}}$ converge. Moreover, $ {\m h}^{\mathrm{R\_TK}} = {\m h}^{\mathrm{set}}$ is a constant vector. Hence, we note that the right hand side is a convergent vector. The size of $[\m A^\mathrm{J}_{\mathrm{P}}\ \m A^\mathrm{J}_{\mathrm{M}}]^\top$ is ${(n_p+n_m)\times n_j}$ and it is clear that the number of pipes and pumps is greater than or equal to the number of junctions in a looped network. That is, the matrix $[\m A^\mathrm{J}_{\mathrm{P}}\ \m A^\mathrm{J}_{\mathrm{M}}]^\top$ has more rows than columns or is a square matrix, and $\m h^\mathrm{J} $ which can be expressed by the convergent vector on the right hand side also converges.
			
			When PRVs or FCVs are in ``ACTIVE" statuses, the proof is exactly the same, because the models of active PRVs and FCVs have similar mathematical form  as the models of tanks or reservoirs and junctions. {For example, when a PRV is active, and the head is set to $\m h^{\mathrm{R}_\mathrm{set}}$. That is, $\m A_{\m h}^\mathrm{W} \m h = \m h^{\mathrm{R}_\mathrm{set}}$ and $\m b^{\mathrm{W}} = \m h^{\mathrm{R}_\mathrm{set}}$ in~\eqref{equ:LinearMatrix}. In fact, this model is exactly the same as the  model of tanks or reservoirs which is $\m A_{\m h}^\mathrm{TK} \m h = \m h^{\mathrm{TK}_\mathrm{set}}$. Similarly, when a FCV is active, and the flow is set to $\m q^{\mathrm{R}_\mathrm{set}}$. That is, $\m A_{\mathrm{W}}^{\mathrm{W}} \m q = \m q^{\mathrm{R}_\mathrm{set}}$ which is similar to mass balance equation $\m A_{\m {q}}^{\mathrm{J}} \m q = \m d$ in~\eqref{equ:LinearMatrix}. It means we can embed the models of active PRVs or FCVs into~\eqref{equ:reogan} directly, and the proof still holds.}
			
		\end{proof}

		

		\normalcolor

		\vspace{-0.75cm}
		
		\section*{Acknowledgments}
		This material is based upon work supported by the National Science Foundation under Grants CMMI-DCSD-1728629 and 1847125.  This work was also supported by the University of Texas at Austin Startup Grant and by Cooperative Agreement No. 83595001 awarded by the U.S. Environmental Protection Agency (EPA) to The University of Texas at Austin. This work has not been formally reviewed by EPA. The views expressed in this document are solely those of the authors and do not necessarily reflect those of the Agency. EPA does not endorse any products or commercial services mentioned in this publication. All the codes, tested networks, and results are freely available on Github for research reproducibility. In order to use GP-NET properly, the reader is referred to the Readme.md file on Github~\cite{shenwangGP}.

		
		%
		%

	\bibliographystyle{IEEEtran}
	\bibliography{IEEEabrv,bibfile.bib}

	\end{document}